\documentclass[preprint,3p,12pt]{elsarticle}
\usepackage{amsmath}
\usepackage{stmaryrd}
\usepackage{bbding}
\usepackage{dcolumn}
\usepackage{graphicx}
\usepackage{amsfonts}
\usepackage{amssymb}
\usepackage{psfrag}
\usepackage{wrapfig}
\usepackage{subfigure}
\usepackage{makeidx}
\usepackage{bm}
\usepackage{epsf}
\usepackage{epsfig}
\usepackage{setspace}
\usepackage{color}
\usepackage{algorithm}
\usepackage{algpseudocode}
\usepackage{setspace}
\usepackage{diagbox}

\begin{document}
\title{Multiple-GPU accelerated high-order gas-kinetic scheme on three-dimensional unstructured meshes}

\author[BNU]{Yuhang Wang}
\ead{hskwyh@outlook.com}

\author[BNU]{Waixiang Cao}
\ead{caowx@bnu.edu.cn}

\author[BNU]{Liang Pan\corref{cor}}
\ead{panliang@bnu.edu.cn}

\address[BNU]{Laboratory of Mathematics and Complex Systems, School of Mathematical Sciences, Beijing Normal University, Beijing, China}
\cortext[cor]{Corresponding author}

\begin{abstract}
Recently, successes have been achieved for the high-order
gas-kinetic schemes (HGKS) on unstructured meshes for compressible
flows. In this paper, to accelerate the computation, HGKS is
implemented with the graphical processing unit (GPU) using the
compute unified device architecture (CUDA).  HGKS on unstructured
meshes is a fully explicit scheme, and the acceleration framework
can be developed based on the cell-level parallelism.  For
single-GPU computation, the connectivity of geometric information is
generated for the requirement of data localization and independence.
Based on such data structure, the kernels and corresponding girds of
CUDA are set. With the one-to-one mapping between the indices of
cells and CUDA threads, the single-GPU computation using CUDA can be
implemented for HGKS. For multiple-GPU computation, the domain
decomposition and data exchange need to be taken into account. The
domain is decomposed into subdomains by METIS, and the MPI processes
are created for the control of each process and communication among
GPUs. With reconstruction of connectivity and adding ghost cells,
the main configuration of CUDA for single-GPU can be inherited by
each GPU. The benchmark cases for compressible flows, including
accuracy test and flow passing through a sphere, are presented to
assess the numerical performance of HGKS with Nvidia RTX A5000 and
Tesla V100 GPUs. For single-GPU computation, compared with the
parallel central processing unit (CPU) code running on the Intel
Xeon Gold 5120 CPU with open multi-processing (OpenMP) directives,
5x speedup is achieved by RTX A5000 and  9x speedup is achieved by
Tesla V100. For multiple-GPU computation, HGKS code scales properly
with the increasing number of GPU. Numerical results confirm the
excellent performance of multiple-GPU accelerated HGKS on
unstructured meshes.
\end{abstract}

\begin{keyword}
High-order gas-kinetic scheme, unstructured meshes, WENO
reconstruction, GPU accelerated computation.
\end{keyword}

\maketitle

\section{Introduction}
Graphical processing unit (GPU) is a form of hardware acceleration,
which is originally developed for graphics manipulation and is
extremely efficient at processing large amounts of data in parallel.
Since these units have a parallel computation capability inherently,
they can provide fast and low cost solutions to high performance
computing (HPC). In recent years, GPUs have gained significant
popularity as a cheaper, more efficient, and more accessible
alternative to large-scale HPC systems with central processing units
(GPU). Great effort has been already achieved for computational
fluid dynamics,  such as direct numerical simulation of turbulent
flows \cite{GPU-1,GPU-2,GPU-3}, high-order CFD simulations with
complex grids \cite{GPU-3}, incompressible smoothed particle
hydrodynamics \cite{GPU-4}, multiphase flows \cite{GPU-5}, shallow
water flows \cite{GPU-6}, unified gas kinetic wave-particle method
\cite{GPU-7} and discrete unified gas kinetic scheme \cite{GPU-8},
etc.

In the past decades, the gas-kinetic scheme (GKS) has been developed
systematically based on the Bhatnagar-Gross-Krook (BGK) model
\cite{BGK-1,BGK-2} under the finite volume framework, and applied
successfully in the computations from low speed flow to hypersonic
one \cite{GKS-Xu1,GKS-Xu2}. The gas-kinetic scheme presents a gas
evolution process from kinetic scale to hydrodynamic scale, where
both inviscid and viscous fluxes are recovered from a time-dependent
and genuinely multi-dimensional gas distribution function at a cell
interface. Starting from a time-dependent flux function, based on
the two-stage fourth-order formulation \cite{GRP-high-1,GRP-high-2},
the high-order gas-kinetic scheme (HGKS) has been constructed and
applied for the compressible flow simulation
\cite{GKS-high-1,GKS-high-2,GKS-high-3}. Originally, the parallel
HGKS code was developed with central processing unit (CPU) using
open multi-processing (OpenMP) directives. However, due to the
limited shared memory, the computational scale is constrained. To
perform the large-scale numerical simulation of turbulence, the
domain decomposition and the message passing interface (MPI)
\cite{MPI-1} are used for parallel implementation \cite{GKS-MPI}. To
further improve the efficiency, the HGKS code is implemented with
single GPU. A major limitation in single-GPU computation is its
available memory, which leads to a bottleneck in the maximum number
of computational mesh. To implement much larger scale computation
and accelerate the efficiency, HGKS is implemented with multiple
GPUs using CUDA and MPI architecture (MPI + CUDA). Due to the
explicit formulation of HGKS, the CPU code with MPI scales properly
with the number of processors used. The numerical results
demonstrates the capability of HGKS as a powerful DNS tool from the
low speed to supersonic turbulence study \cite{GKS-MPI}. Recently,
the three-dimension discontinuous Galerkin based HGKS has been
implemented in single-GPU computation using compute unified device
architecture (CUDA) \cite{GKS-GPU-1,GKS-GPU-2}. Obtained results are
compared with those obtained by Intel i7-9700 CPU using OpenMP
directives. The GPU code achieves 6x-7x speedup with TITAN RTX, and
10x-11x speedup with Tesla V100. The computational time of parallel
CPU code running on $1024$ Intel Xeon E5-2692 cores with MPI is
approximately $3$ times longer than that of GPU code using $8$ Tesla
V100 GPUs with MPI and CUDA.

To deal with the complicated geometry, the high-order gas-kinetic
scheme on unstructured meshes has been developed with the WENO
reconstruction \cite{GKS-high-3,GKS-high-4}. Due to the complex mesh
topology of hybrid unstructured meshes, the procedures of classical
WENO scheme, including the selection of candidate stencils and
calculation of linear weights, become extremely complicated. A
simple strategy of selecting stencils for reconstruction is adopted
and the topology independent linear weights are used. A large
stencil is selected with the neighboring cells and the neighboring
cells of neighboring cells, and a quadratic polynomial can be
obtained.  A robust selections of candidate sub-stencils are also
given, such that the linear polynomials are solvable for each
sub-stencil. In this paper, to accelerate the computation, HGKS  on
unstructured meshes is implemented with GPU using CUDA. The current
HGKS is a fully explicit scheme,  and the acceleration of
computation can be achieved by executing the calculation of the
cells simultaneously. Compared with the structured meshes,  the
connectivity of geometric information for unstructured meshes needs
to be generated for the requirement of data localization and
independence. Based on such data structure, the kernels and
corresponding girds of CUDA are set. With the one-to-one mapping
between the indices of cells and CUDA threads, the single-GPU
computation using CUDA can be implemented for HGKS. For multiple-GPU
computation, the domain decomposition and data exchange, which are
more complicated than that of structured meshes, need to be taken
into account.  The domain is decomposed into several subdomains by
METIS \cite{METIS-Kar}, and the MPI processes are created for
one-to-one GPU management and communication. The CUDA-Aware MPI
library \cite{CUDA-MPI} is used for GPU-GPU communication. With
reconstruction of connectivity and adding three-layer of ghost
cells, the main configuration of CUDA for single-GPU can be
inherited by each GPU. The mappings between two subdomains are built
for the transfer of conservative variables of ghost cells.
Non-blocking communication and the corresponding strategy are used
to improve the efficiency of GPU-GPU communication. For single-GPU
implementation using CUDA, compared with the CPU code using 2
Intel(R) Xeon(R) Gold 5120 CPUs with OpenMP directives, 5x speedup
is achieved for RTX A5000 and 9x speedup is achieved for Tesla V100.
For multiple-GPU with CUDA and MPI, the HGKS is strongly scalable
with the increasing number of GPUs.  Nearly linear speedup can be
achieved under suitable computational work-load. Numerical
performance shows that the data communication crossing GPUs through
MPI costs the relative little time. To reduce the memory cost and
improve the computational efficiency,  the code of multiple-GPU
accelerated HGKS is compiled using FP32 (single) precision for the
accuracy test. As expected, the efficiency can be improved and the
memory cost can be reduced with FP32 precision. Compared with the
results of FP64 (double) precision,  the errors of the accuracy
increase slightly and the third-order can be maintained.

This paper is organized as follows. In Section 2, the high-order
gas-kinetic scheme is briefly reviewed. The single-GPU
implementation for HGKS code is introduced in Section 3. Section 4
includes the multiple-GPU implementation. The numerical results are
presented in Section 5. The last section is the conclusion.

\section{High-order gas-kinetic scheme}
\subsection{BGK equation and finite volume scheme}
The Boltzmann equation expresses the behavior of a many-particle
kinetic system in terms of the evolution equation for a single
particle gas distribution function. The BGK equation
\cite{BGK-1,BGK-2}  is the simplification of Boltzmann equation, and
the three-dimensional BGK equation can be written as
\begin{equation}\label{bgk}
f_t+uf_x+vf_y+wf_z=\frac{g-f}{\tau},
\end{equation}
where $\boldsymbol{u}=(u,v,w)$ is the particle velocity, $\tau$ is
the collision time, $f$ is the gas distribution function. $g$ is the
equilibrium state given by Maxwellian distribution
\begin{equation*}
g=\rho(\frac{\lambda}{\pi})^{(N+3)/2}e^{-\lambda[(u-U)^2+(v-V)^2+(w-W)^2+\xi^2]},
\end{equation*}
where $\rho$ is the density, $\boldsymbol{U}=(U,V,W)$ is the
macroscopic fluid velocity, and $\lambda$ is the inverse of gas
temperature, i.e., $\lambda=m/2kT$. In the BGK model, the collision
operator involves a simple relaxation from $f$ to the local
equilibrium state $g$. The variable $\xi$ accounts for the internal
degree of freedom, $\xi^2=\xi_1^2+\dots+\xi_N^2$,
$N=(5-3\gamma)/(\gamma-1)$ is the internal degree of freedom, and
$\gamma$ is the specific heat ratio. The collision term satisfies
the compatibility condition
\begin{equation*}
 \int \frac{g-f}{\tau}\psi \text{d}\Xi=0,
\end{equation*}
where $\displaystyle\psi=(1,u,v,w,\frac{1}{2}(u^2+v^2+w^2+\xi^2))^T$
and
$\text{d}\Xi=\text{d}u\text{d}v\text{d}w\text{d}\xi_1\dots\text{d}\xi_{N}$.
According to the Chapman-Enskog expansion for BGK equation, the
macroscopic governing equations can be derived. In the continuum
region, the BGK equation can be rearranged and the gas distribution
function can be expanded as
\begin{align*}
f=g-\tau D_{\boldsymbol{u}}g+\tau D_{\boldsymbol{u}}(\tau
D_{\boldsymbol{u}})g-\tau D_{\boldsymbol{u}}[\tau
D_{\boldsymbol{u}}(\tau D_{\boldsymbol{u}})g]+...,
\end{align*}
where $D_{\boldsymbol{u}}=\displaystyle\frac{\partial}{\partial
t}+\boldsymbol{u}\cdot \nabla$. With the zeroth-order truncation
$f=g$, the Euler equations can be obtained. For the first-order
truncation
\begin{align*}
f=g-\tau (ug_x+vg_y+wg_z+g_t),
\end{align*}
the Navier-Stokes equations can be obtained \cite{GKS-Xu1,GKS-Xu2}.

Taking moments of Eq.\eqref{bgk} and integrating with respect to
space, the semi-discretized finite volume scheme can be expressed as
\begin{align}\label{semi}
\frac{\text{d} Q_i}{\text{d} t}=\mathcal{L}(Q_{i}),
\end{align}
where $Q_i=(\rho, \rho U,\rho V, \rho W, \rho E)$ is the cell
averaged conservative value of $\Omega_{i}$, $\rho$ is the density,
$U,V,W$ is the flow velocity and $\rho E$ is the total energy
density. The operator $\mathcal{L}$ is defined as
\begin{equation}\label{finite}
\mathcal{L}(Q_{i})=-\frac{1}{|\Omega_{i}|}\sum_{i_p\in
N(i)}F_{i,i_p}(t)S_{i_p}=-\frac{1}{|\Omega_{i}|}\sum_{i_p\in
N(i)}\iint_{\Sigma_{i_p}}\boldsymbol{F}(Q,t)\text{d}\sigma,
\end{equation}
where $|\Omega_{i}|$ is the volume of $\Omega_{i}$, $\Sigma_{i_p}$
is the common cell interface of $\Omega_{i}$, $S_{i_p}$ is the area
of $\Sigma_{i_p}$ and $N(i)$ is the set of index for neighboring
cells of $\Omega_{i}$. To achieve the expected order of accuracy,
the Gaussian quadrature is used for the flux integration
\begin{align*}
\iint_{\Sigma_{i_p}}\boldsymbol{F}(Q,t)\text{d}\sigma=\sum_{G}\omega_{G}F_{G}(t)S_{i_p},
\end{align*}
where $\omega_{G}$ is the quadrature weights. The numerical flux
$F_{G}(t)$ at Gaussian quadrature point $\boldsymbol{x}_{G}$ can be
given by taking moments of gas distribution function
\begin{align}\label{flux-G}
F_{G}(t)=\int\boldsymbol\psi \boldsymbol{u}\cdot\boldsymbol{n}_{G}
f(\boldsymbol{x}_{G},t,\boldsymbol{u},\xi)\text{d}\Xi,
\end{align}
where $F_{G}(t)=(F^{\rho}_{G},F^{\rho U}_{G},F^{\rho V}_{G},F^{\rho
W}_{G}, F^{\rho E}_{G})$ and $\boldsymbol{n}_{G}$ is the local
normal direction of  cell interface. With the coordinate
transformation, the numerical flux in the global coordinate can be
obtained. Based on the integral solution of BGK equation
Eq.\eqref{bgk}, the gas distribution function
$f(\boldsymbol{x}_{G},t,\boldsymbol{u},\xi)$ in the local coordinate
can be given by
\begin{equation*}
f(\boldsymbol{x}_{G},t,\boldsymbol{u},\xi)=\frac{1}{\tau}\int_0^t
g(\boldsymbol{x}',t',\boldsymbol{u},
\xi)e^{-(t-t')/\tau}\text{d}t'+e^{-t/\tau}f_0(-\boldsymbol{u}t,\xi),
\end{equation*}
where $\boldsymbol{x}'=\boldsymbol{x}_{G}-\boldsymbol{u}(t-t')$ is
the trajectory of particles, $f_0$ is the initial gas distribution
function, and $g$ is the corresponding equilibrium state. With the
first order spatial derivatives, the second-order gas distribution
function at cell interface can be expressed as
\begin{equation}\label{flux}
\begin{split}
f(\boldsymbol{x}_{G},t,\boldsymbol{u},\xi)= & (1-e^{-t/\tau})g_0
+((t+\tau)e^{-t/\tau}-\tau)(\overline{a}_1u+\overline{a}_2v+\overline{a}_3w)g_0\\
+& (t-\tau+\tau e^{-t/\tau}){\bar{A}} g_0 \\
+& e^{-t/\tau}g_r[1-(\tau+t)(a_{1}^{r}u+a_{2}^{r}v+a_{3}^{r}w)-\tau A^r)](1-H(u)) \\
+& e^{-t/\tau}g_l[1-(\tau+t)(a_{1}^{l}u+a_{2}^{l}v+a_{3}^{l}w)-\tau A^l)]H(u),
\end{split}
\end{equation}
where the equilibrium state $g_{0}$ and the corresponding
conservative variables $Q_{0}$ can be determined by the
compatibility condition
\begin{align*}
\int\psi g_{0}\text{d}\Xi=Q_0=\int_{u>0}\psi
g_{l}\text{d}\Xi+\int_{u<0}\psi g_{r}\text{d}\Xi.
\end{align*}
With the reconstruction of macroscopic variables, the coefficients
in Eq.\eqref{flux} can be fully determined by the reconstructed
derivatives and compatibility condition
\begin{equation*}
    \begin{aligned}
        \displaystyle
        \langle a_{1}^{k}\rangle=\frac{\partial Q_{k}}{\partial \boldsymbol{n_x}},
        \langle a_{2}^{k}\rangle=\frac{\partial Q_{k}}{\partial \boldsymbol{n_y}},
        \langle a_{3}^{k}\rangle     & =\frac{\partial Q_{k}}{\partial\boldsymbol{n_z}},
        \langle
        a_{1}^{k}u+a_{2}^{k}v+a_{3}^{k}w+A^{k}\rangle=0,                                  \\
        \displaystyle
        \langle\overline{a}_1\rangle=\frac{\partial Q_{0}}{\partial \boldsymbol{n_x}},
        \langle\overline{a}_2\rangle=\frac{\partial Q_{0}}{\partial \boldsymbol{n_y}},
        \langle\overline{a}_3\rangle & =\frac{\partial Q_{0}}{\partial \boldsymbol{n_z}},
        \langle\overline{a}_1u+\overline{a}_2v+\overline{a}_3w+\overline{A}\rangle=0,
    \end{aligned}
\end{equation*}
where $k=l$ and $r$,  $\boldsymbol{n_x}$, $\boldsymbol{n_y}$,
$\boldsymbol{n_z}$ are the unit directions of local coordinate at
$\boldsymbol{x}_{G}$ and $\langle...\rangle$ are the moments of the
equilibrium $g$ and defined by
\begin{align*}
\langle...\rangle=\int g (...)\psi \text{d}\Xi.
\end{align*}
More details of the gas-kinetic scheme can be found in
\cite{GKS-Xu1,GKS-Xu2}.

\subsection{Temporal reconstruction}
In this paper,  the two-stage fourth-order temporal discretization
\cite{GRP-high-1,GKS-high-1} is used to achieve the high-order
temporal accuracy. For the semi-discretized finite volume scheme
Eq.\eqref{semi}, the flow variables $Q^{n+1}$ at $t_{n+1}=t_n+\Delta
t$  can be updated with the following formula
\begin{align*}
& Q^*_i=Q^n_i+\frac{1}{2}\Delta
t\mathcal{L}(Q^n_i)+\frac{1}{8}\Delta
    t^2\frac{\partial}{\partial
    t}\mathcal{L}(Q^n_i),                                                             \\
    Q^{n+1}_{i} & =Q^n_i+\Delta t\mathcal {L}(Q^n_i)+\frac{1}{6}\Delta
    t^2\big(\frac{\partial}{\partial
            t}\mathcal{L}(Q^n_i)+2\frac{\partial}{\partial
            t}\mathcal{L}(Q^*_i)\big).
\end{align*}
It can be proved that the above temporal discretization provides a
fourth-order time accurate solution for hyperbolic equations
\cite{GRP-high-1}. To implement the two-stage method for
Eq.\eqref{semi}, a linear function is used to approximate the time
dependent numerical fluxes  Eq.\eqref{flux-G}
\begin{align}\label{expansion-1}
    F_{G}(t)=F_{G}^n+ \partial_t F_{G}^n(t-t_n).
\end{align}
Integrating Eq.\eqref{expansion-1} over $[t_n, t_n+\Delta t/2]$ and
$[t_n, t_n+\Delta t]$, we have the following two equations
\begin{align*}
    F_{G}^n\Delta t            & +\frac{1}{2}\partial_t F_{G}^n\Delta t^2 =\int_{t_n}^{t_n+\Delta t}F_{G}(t)\text{d}t,   \\
    \frac{1}{2}F_{G}^n\Delta t & +\frac{1}{8}\partial_t F_{G}^n\Delta t^2 =\int_{t_n}^{t_n+\Delta t/2}F_{G}(t)\text{d}t.
\end{align*}
The coefficients $F_{G}^n$ and $\partial_t F_{G}^n$ can be
determined by solving the linear system above, and the operator
$\mathcal{L}$ and its temporal derivative $\partial_t\mathcal{L}$ at
$t=t_n$ can be given according to Eq.\eqref{finite}. Similarly,
$\mathcal{L}$ and $\partial_t\mathcal{L}$  at the intermediate state
$t=t_n+\Delta t/2$ can be constructed as well.

\subsection{Spatial reconstruction}
To deal with the complex geometry, the three-dimensional
unstructured meshes are considered, and the tetrahedral and
hexahedral meshes are used in this paper for simplicity. In the
previous studies, the high-order gas-kinetic schemes have be
developed with the third-order non-compact WENO reconstruction
\cite{GKS-high-3, GKS-high-4,WENO-simple-1,WENO-simple-2}. In this
section, a brief review of  WENO reconstruction on unstructured
hexahedral and tetrahedral meshes will be presented.

For the cell $\Omega_i$, the faces are labeled as $f_{p}$, where
$p=1,\dots,4$ for tetrahedral cell, and $p=1,\dots,6$ for hexahedral
cell. The neighboring cell of $\Omega_{i}$, which shares the face
$f_{p}$, is denoted as $\Omega_{i_p}$. Meanwhile, the neighboring
cells of $\Omega_{i_p}$ are denoted as $\Omega_{i_{pm}}$. For WENO
reconstruction, a big stencil  for cell $\Omega_i$  is selected as
follows
\begin{align*}
S_i^{WENO} & =\{\Omega_{i},\Omega_{i_p},\Omega_{i_{pm}}\},
\end{align*}
which is consist of neighboring cells and neighboring cells of
neighboring cells of $\Omega_{i}$. To deal with the discontinuity,
the sub-stencils $S_{i_m}^{WENO}$ in non-compact WENO reconstruction
for cell $\Omega_{i}$ are selected, where $m=1,\dots, M$ and $M$ is
the number of sub-stencils. For the hexahedral cell, $M=8$ and the
sub-candidate stencils are selected as
\begin{align*}
    S_{i_1}^{WENO}=\{\Omega_{i},\Omega_{i_1},\Omega_{i_2},\Omega_{i_3}\},~S_{i_5}^{WENO}=\{\Omega_{i},\Omega_{i_6},\Omega_{i_2},\Omega_{i_3}\}, \\
    S_{i_2}^{WENO}=\{\Omega_{i},\Omega_{i_1},\Omega_{i_3},\Omega_{i_4}\},~S_{i_6}^{WENO}=\{\Omega_{i},\Omega_{i_6},\Omega_{i_3},\Omega_{i_4}\}, \\
    S_{i_3}^{WENO}=\{\Omega_{i},\Omega_{i_1},\Omega_{i_4},\Omega_{i_5}\},~S_{i_7}^{WENO}=\{\Omega_{i},\Omega_{i_6},\Omega_{i_4},\Omega_{i_5}\}, \\
    S_{i_4}^{WENO}=\{\Omega_{i},\Omega_{i_1},\Omega_{i_5},\Omega_{i_2}\},~S_{i_8}^{WENO}=\{\Omega_{i},\Omega_{i_6},\Omega_{i_5},\Omega_{i_2}\}.
\end{align*}
The linear polynomials can be determined based on above stencils,
which contain the target cell $\Omega_{i}$ and three neighboring
cells. For the tetrahedral cells, in order to avoid the centroids of
$\Omega_{i}$ and three of neighboring cells becoming coplanar,
additional cells are needed for the sub-candidate stencils. For the
tetrahedral cell,  four sub-candidate stencils are selected as
\begin{align*}
    S_{i_1}^{WENO}=\{\Omega_{i},\Omega_{i_1},\Omega_{i_2},\Omega_{i_3},\Omega_{i_{11}},\Omega_{i_{12}},\Omega_{i_{13}}\}, \\
    S_{i_2}^{WENO}=\{\Omega_{i},\Omega_{i_1},\Omega_{i_2},\Omega_{i_4},\Omega_{i_{21}},\Omega_{i_{22}},\Omega_{i_{23}}\}, \\
    S_{i_3}^{WENO}=\{\Omega_{i},\Omega_{i_2},\Omega_{i_3},\Omega_{i_4},\Omega_{i_{31}},\Omega_{i_{32}},\Omega_{i_{33}}\}, \\
    S_{i_4}^{WENO}=\{\Omega_{i},\Omega_{i_3},\Omega_{i_1},\Omega_{i_4},\Omega_{i_{41}},\Omega_{i_{42}},\Omega_{i_{43}}\}.
\end{align*}
The cells of sub-candidate stencils are consist of the three
neighboring cells and three neighboring cells of one neighboring
cell. With such an enlarged sub-stencils, the linear polynomials can
be determined.

With the selected stencil, a quadratic polynomial and several linear
polynomials can be constructed based on the big stencil and the
sub-stencils as follows
\begin{equation}\label{polys}
    \begin{split}
        P_0(\boldsymbol{x})&=Q_{0}+\sum_{|\boldsymbol d|=1}^2a_{\boldsymbol d}p_{\boldsymbol d}(\boldsymbol{x}),\\
        P_m(\boldsymbol{x})&=Q_{0}+\sum_{|\boldsymbol d|=1}b_{\boldsymbol d}^mp_{\boldsymbol d}(\boldsymbol{x}),
    \end{split}
\end{equation}
where $m=1,\dots,M$, $Q_{0}$ is the cell averaged variables over
$\Omega_{0}$ with newly rearranged index, the multi-index
$\boldsymbol d=(d_1, d_2,d_3)$ and $|\boldsymbol d|=d_1+d_2+d_3$.
The base function $p_{\boldsymbol d}(\boldsymbol{x})$ is defined as
\begin{align*}
\displaystyle p_{\boldsymbol
d}(\boldsymbol{x})=x^{d_1}y^{d_2}z^{d_3}-\frac{1}{\left|\Omega_{0}\right|}\iiint_{\Omega_{0}}x^{d_1}y^{d_2}z^{d_3}\text{d}V.
\end{align*}
To determine these polynomials for WENO reconstruction, the following constrains need to be satisfied
\begin{align*}
\frac{1}{\left|\Omega_{k}\right|}\iiint_{\Omega_{k}}P_0(\boldsymbol{x})\text{d}V   & =Q_{k},~\Omega_{k}\in S_i^{WENO},     \\
\frac{1}{\left|\Omega_{m_k}\right|}\iiint_{\Omega_{k}}P_m(\boldsymbol{x})\text{d}V
& =Q_{k},~\Omega_{k}\in S_{i_m}^{WENO},
\end{align*}
where $k=0, \dots, K$, $m_k=0,\dots,m_K$, $Q_{k}$ and
$Q_{m_k}$ are the conservative variable with newly rearranged index.
The over-determined linear systems can be generated and the least
square method is used to obtain the coefficients $a_{\boldsymbol d}$
and $b_{\boldsymbol d}^m$.

With the reconstructed polynomial $P_m(\boldsymbol{x}), m=0,...,M$,
the point-value $Q(\boldsymbol{x}_{G})$ and the spatial derivatives
$\partial_{x,y,z} Q(\boldsymbol{x}_{G})$  for reconstructed
variables  at Gaussian quadrature point can be given by the
non-linear combination
\begin{equation}\label{weno}
    \begin{split}
        Q(\boldsymbol{x}_{G})=\overline{\omega}_0(\frac{1}{\gamma_0}P_0(\boldsymbol{x}_{G})-&
        \sum_{m=1}^{M}\frac{\gamma_m}{\gamma_0}P_m(\boldsymbol{x}_{G}))+\sum_{m=1}^{M}\overline{\omega}_mP_m(\boldsymbol{x}_{G}),\\
        \partial_{x,y,z} Q(\boldsymbol{x}_{G})=\overline{\omega}_0(\frac{1}{\gamma_0}\partial_{x,y,z}
        P_0(\boldsymbol{x}_{G})-&\sum_{m=1}^{M}\frac{\gamma_m}{\gamma_0}\partial_{x,y,z}
        P_m(\boldsymbol{x}_{G}))+\sum_{m=1}^{M}\overline{\omega}_m\partial_{x,y,z}
        P_m(\boldsymbol{x}_{G}),
    \end{split}
\end{equation}
where $\gamma_0, \gamma_1,\dots,\gamma_M$ are the linear weights.
The non-linear weights $\omega_m$ and normalized non-linear weights
$\overline{\omega}_m$ are defined as
\begin{align*}
    \overline{\omega}_{m}=\frac{\omega_{m}}{\sum_{m=0}^{M} \omega_{m}},~
    \omega_{m}=\gamma_{m}\Big[1+\Big(\frac{\tau_Z}{\beta_{m}+\epsilon}\Big)\Big],~
    \tau_Z=\sum_{m=1}^{M}\Big(\frac{|\beta_0-\beta_m|}{M}\Big),
\end{align*}
where $\epsilon$ is a small positive number. The smooth indicator
$\beta_{m}$ is defined as
\begin{align*}
    \beta_m=\sum_{|l|=1}^{r_m}|\Omega_{i}|^{\frac{2|l|}{3}-1}\int_{\Omega_{i}}
    \Big(\frac{\partial^lP_m}{\partial_x^{l_1}\partial_y^{l_2}\partial_z^{l_3}}(x,y,z)\Big)^2\text{d}V,
\end{align*}
where $r_0=2$ and $r_m=1$ for $m=1,\dots,M$.  It can be proved that
Eq.\eqref{weno} ensures third-order accuracy and more details can be
found in \cite{GKS-high-3, GKS-high-4}. In the computation, the
linear weights are set as $\gamma_i=0.025$, $\gamma_0=1-\gamma_i M$
for both non-compact and compact scheme without special statement.

\section{Single-GPU implementation}
With the development of hardware and software, GPUs are becoming an
important part of high performance computing, which provide great
computational power for large-scale scientific problem.  In the
previous studies \cite{GKS-GPU-1,GKS-GPU-2}, the parallel
acceleration for HGKS on structured meshes is implemented by GPUs
with Compute Unified Device Architecture (CUDA) and Message Passing
Interface (MPI).  The multiple-GPU accelerated HGKS code scales
properly with the increasing number of GPU. In this paper, the
multiple-GPU accelerated HGKS code will be extended to unstructured
meshes. HGKS on unstructured meshes is a fully explicit scheme, and
the idea of parallelism at the cell level is retained, which means
that  the acceleration can be achieved by performing calculations on
a large number of cells simultaneously. In this section, we mainly
focus on the single-GPU implementation with CUDA. The multiple-GPU
implementation with MPI and CUDA with be presented in next section.

\begin{figure}[!h]
\centering
\includegraphics[width=0.5\textwidth]{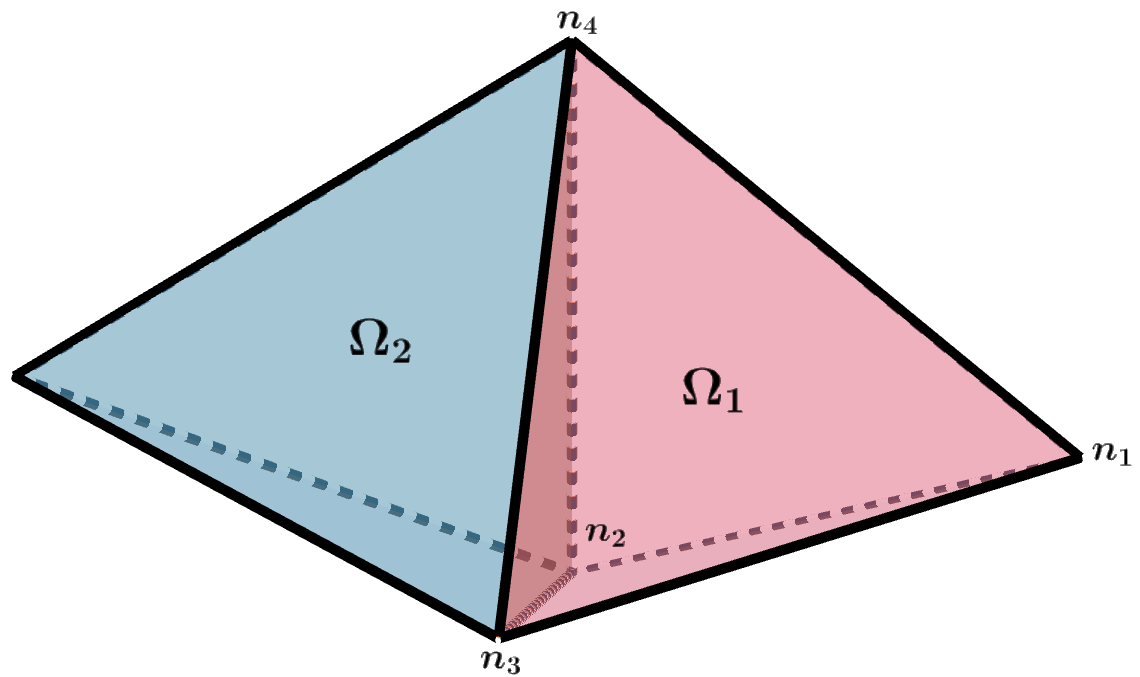}
\caption{\label{uns-single-geometric-1} Data structure of unstructured meshes in HGKS.}
\end{figure}

\subsection{Generation of geometric information}
For the structured meshes,  the computational domain is consist of
cells, which are distinguished by the indices $(i,j,k)$. These
indices are directly connected with the way of storage for flow
field and geometric data in the computer. The topological
relationships for computational mesh (node, face and cell) can be
obtained naturally. Thus, the data localization and independence are
easily achieved. However, it is not natural for the unstructured
meshes, in which all geometric information is stored in
one-dimensional arrays. Each cell, face, node owns a unique index,
and the adjacent indices cannot be provided directly. Therefore, the
additional geometric connectivities need to be built for the
cell-level parallelization of HGKS.

Three main parts of algorithm are taken into consideration: WENO
reconstruction, GKS flux solver, temporal discretization. To update
the conservative variables of a control volume, the numerical fluxes
and its temporal derivatives at the cell interfaces need be
provided. To calculate the fluxes and its temporal derivatives, the
conservative variables of two-layer neighbors should be accessed for
WENO reconstruction. Besides, the node coordinates for each cell are
needed for the generation of geometric information, such as the
area, volume and local coordinate systems. Thus, the cell-to-node,
face-to-node, cell-to-face, face-to-cell, and and cell-to-neighbor
connectivities need to be constructed. As shown in
Figure.\ref{uns-single-geometric-1}, the tetrahedral cell is taken
as an example, and these geometric connectivities can be built as
follows
\begin{enumerate}
\item Cell-to-node and face-to-node connectivities can be provided by two
two-dimensional arraies, i.e.,  \textit{FaceNode} and
\textit{CellNode}, respectively. As shown in
Figure.\ref{uns-single-geometric-1}, there are four nodes for cell
$\Omega_1$, i.e., $n_1, n_2, n_3, n_4$, and four faces, i.e.,
$f_{1}=\{n_2,n_3,n_4\}, f_{2}=\{n_1,n_3,n_4\},
f_{3}=\{n_1,n_2,n_4\},f_{4}=\{n_1,n_2,n_3\}.$ Therefore, the
connectivities are described by \textit{CellNode}[$\Omega_1$][$4$] =
\{$n_1, n_2, n_3, n_4$\} and  \textit{FaceNode}[$f_1$][$3$] =
\{$n_2,n_3,n_4$\}, where the second index of two arrays should be
changes according to the type of cell.
\item Cell-to-face and  face-to-cell connectivities is defined by two
two-dimensional arraies, i.e., \textit{CellFace} and
\textit{FaceCell}. As shown in Figure.\ref{uns-single-geometric-1},
the cell $\Omega_1$ is composed by four faces $(f_1, f_2, f_3, f_4)$
and the cell $\Omega_1$ and $\Omega_2$ share the common face $f_1$.
The connectivity is described  by \textit{CellFace}[$\Omega_1$][$4$]
= \{$f_1, f_2, f_3, f_4$\} and \textit{FaceNeighbor}[$f_1$][2] =
\{$\Omega_1$,$\Omega_2$\}. The second index of \textit{CellFace}
varies according to the type of cell. For the boundary face $f^*$,
\textit{FaceNeighbor}[$f^*$][2] = \{$\Omega_*$,$0$\}, where $0$ is a
label for boundary cell.
\item Cell-to-neighbor connectivity can be also defined by a
two-dimensional array, i.e., \textit{CellNeighbor}.Compared with the
connectivities above, it can be built by a traversal search with the
cell-to-face and face-to-cell connectivities. \textit{CellFace}
provides the index of faces belonging to each cell, and the
neighbors can be obtained by traverse corresponding
\textit{FaceCell}. According to the data stored in
\textit{FaceCell},  $0$ will be received for the boundary cell and
it will be replaced by the cell rank obtained by using boundary
condition. With the procedure above, the cell-to-neighbor
connectivity is obtained by
\textit{CellNeighbor}[$\Omega$][$\Omega_1, \Omega_2, …,
\Omega_n$], where $n$ varies according to the type of meshes.
\item With the cell-to-neighbor connectivity of neighboring cells of
$\Omega$, the collection of two-layer neighbors of $\Omega$ can be
obtained. Based on the collection of two-layer neighbors, HGKS on
unstructured meshes can be implemented with WENO reconstruction. The
process of building such connectivities is shown in
Algorithm.\ref{alg::two-layer build}.
\end{enumerate}

\begin{algorithm}
        \caption{Build relationship between two-layer neighbors and cells}
        \label{alg::two-layer build}
        \begin{algorithmic}
        \For{CellRank = $1$ to NumTotalCell}
            \For{$\text{FirstLayerCellRank} = \text{1 to NumTotalCellNeibor}$}
                \State $\text{FNeighborID} \leftarrow \textit{CellNeighbor}[\text{CellRank}][\text{FirstLayerCellRank}]$
                \For{$\text{SecondLayerCellRank} = \text{1 to NumTotalFirstNeibor}$}
                \State $\text{SNeighborID} \leftarrow \textit{CellNeighbor}[\text{FNeighborID}][\text{SecondLayerCellRank}]$
                \If{SNeighborID $== 0$}
                    \State Add ghost cell according to boundary condition
                \EndIf
                \If{SNeighborID does not exist in \textit{CellNeighbor}[CellRank][*]}
                    \State Add SNeighborID $\rightarrow$ \textit{CellNeighbor}[\text{CellRank}][*]
                \EndIf
                \EndFor
            \EndFor
        \EndFor
        \end{algorithmic}
    \end{algorithm}

\begin{figure}[!h]
    \centering
    \includegraphics[width=0.6\textwidth]{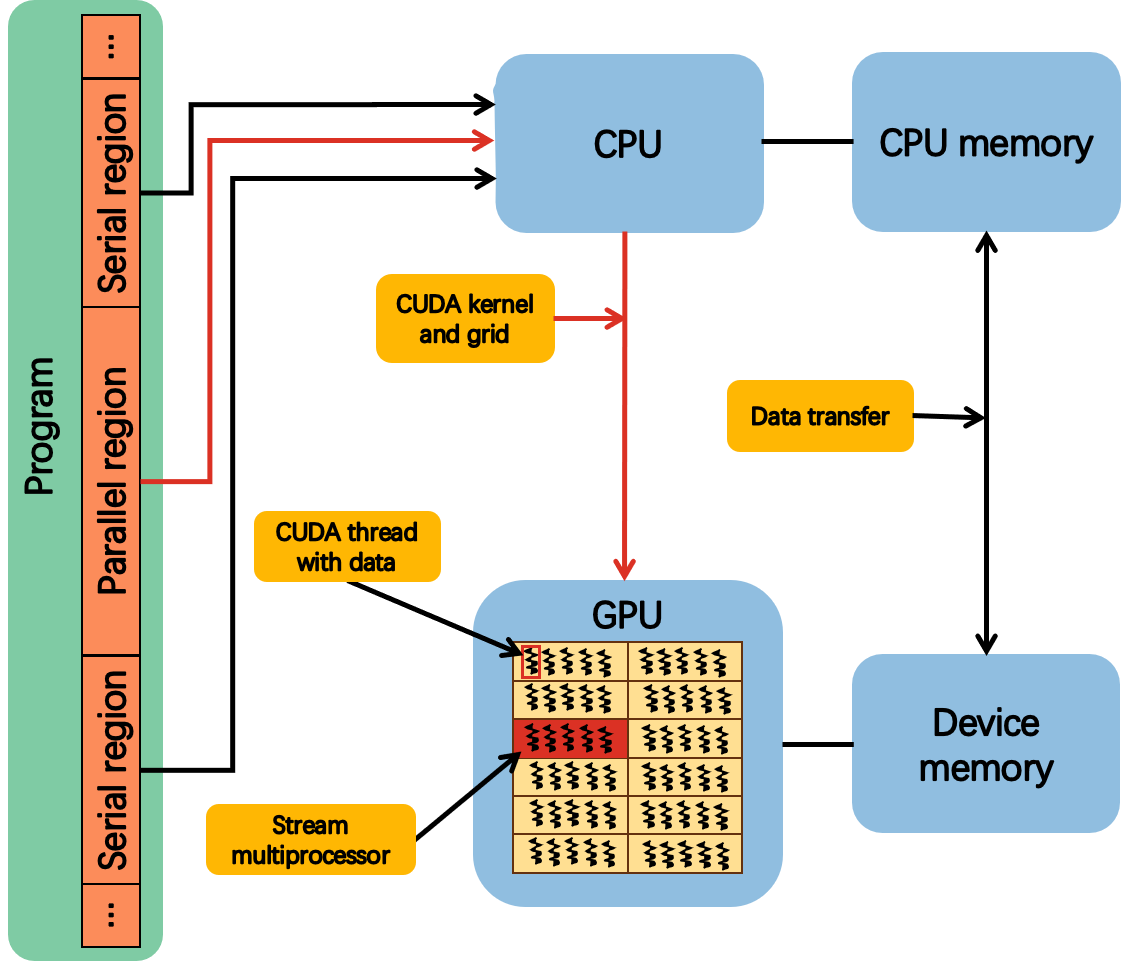}
    \caption{\label{uns-single-CUDA-1}Framework of CUDA with GPU.}
\end{figure}

\subsection{CUDA implementation}
With the connectivities above, the cell-level parallelism will be
implemented with single GPU and Compute unified device architecture
(CUDA). CUDA is a parallel computing platform developed by NVIDIA
for general-purpose computing with GPU. It provides the parallel
computing architecture that introduces a novel programming model
established on high abstraction levels. From the perspective of
CUDA, GPU is viewed as a computer device capable of executing
thousands of threads simultaneously and it works as  a co-processor
to the main CPU. The brief framework of such process is shown in
Figure.\ref{uns-single-CUDA-1}. The CPU is regarded as host, and GPU
is treated as device. Data-parallel, compute-intensive operations
running on the host are transferred to device by using kernels, and
kernels are executed on the device by many different threads. For
CUDA, these threads are organized into thread blocks, and thread
blocks constitute a grid. Such computational structures build
connection with Nvidia GPU hardware architecture. The Nvidia GPU
consists of multiple streaming multiprocessors (SMs), and each SM
contains streaming processors (SPs). When invoking a kernel, the
blocks of grid are distributed to SMs, and the threads of each block
are executed by SPs. The key of parallel computation using CUDA and
single GPU is setting kernels and corresponding grids.

\begin{algorithm}
    \caption{HGKS on unstructured meshes}
    \label{alg:hgk-uns}
    \begin{algorithmic}
        \State Load unstructured meshes information
        \State Establish geometric connectivity in data structure
        \State Build two-layer neighbors with boundary conditions
        \While{Time $\leqslant$ StopTime}
            \State Calculate time step
            \For{Step $= 1$ to $2$}
                \State WENO reconstruction
                \State Compute flux
                \State Update flow variables
            \EndFor
            \If{Time $==$ OutputTime}
                \State Export of output files
            \EndIf
        \EndWhile
    \end{algorithmic}
\end{algorithm}

\begin{figure}[!h]
    \centering
    \includegraphics[width=0.8\textwidth]{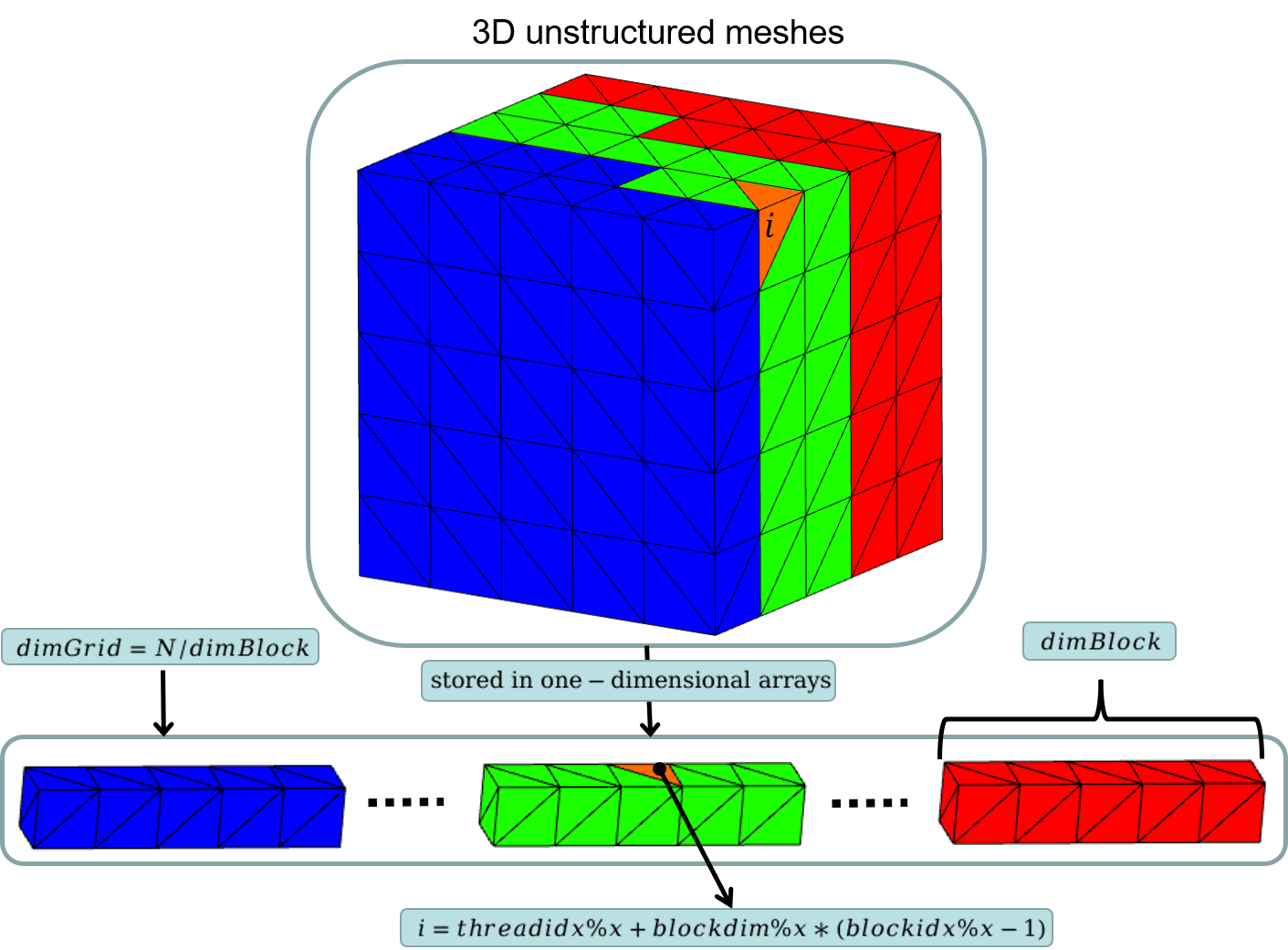}
    \caption{\label{uns-single-grid} Process of setting grid for unstructured meshes.}
\end{figure}

For the explicit computation, HGKS on unstructured meshes is fully
parallelized, and the main process of HGKS code is introduced in
Algorithm.\ref{alg:hgk-uns}. Except for the data transfer involved
in "Load unstructured meshes information" and "Export of output
files",  the rest parts are parallel regions and can be set as
kernels. To achieve  the cell-level parallelism of the algorithm,
the configuration of grids requires the equivalence of the number of
threads and number of cells. The indices of threads should be mapped
one-to-one with the indices of cells.  The kernel invocation can be
described as
\begin{equation*}
kernel<<< dimGrid,dimBlock >>>,
\end{equation*}
where $dimGrid$ and $dimBlock$ are either integer expressions (for
one-dimensional grids and thread blocks), or integer arrays (for
two-dimensional or three-dimensional grids and thread blocks).
$dimGrid$ describes the number of thread blocks in the grid, and the
number of threads in each thread block can be determined by
$dimBlock$. The brief process of setting grids is shown in
Figure.\ref{uns-single-grid}. For unstructured meshes, the indices
of cells are provided by a one-dimensional array, and $dimGrid$ and
$dimBlock$ will be set as two integers. A single thread block, which
contains the same number of threads as the cells, should be used in
the grid. However, there are some limitations on the number of
threads in each thread block, i.e., each thread block contains a
maximum of 1024 threads. Assume the total number of cells is $N$,
$dimGrid$ is defined as $dimGrid = N/dimBlock$.  If $N$ is not
divisible by $dimBlock$, an extra thread block is needed. The
relationship between  the cells indices $i$ and threads indices can
be defined as
\begin{equation*}
    i=threadidx\%x+blockdim\%x*(blockidx\%x-1).
\end{equation*}
Thus, HGKS code can be implemented with specifying kernels and grids with single GPU.

\section{Multiple-GPU implementation}
For the HGKS code on structured meshes \cite{GKS-GPU-2}, the
one-dimensional decomposition is adopted, and MPI processes are
created to manage GPU computation and data transfer.  The number of
parts equals to the number of GPUs. Each GPU inherits the
configuration of single GPU implementation to complete parallel
computation. In this paper, such framework will be extended to HGKS
on unstructured meshes.  The implementation is presented in this
section, including domain decomposition and the combination of MPI
and CUDA.

\begin{figure}[!h]
\centering
\includegraphics[width=0.9\textwidth]{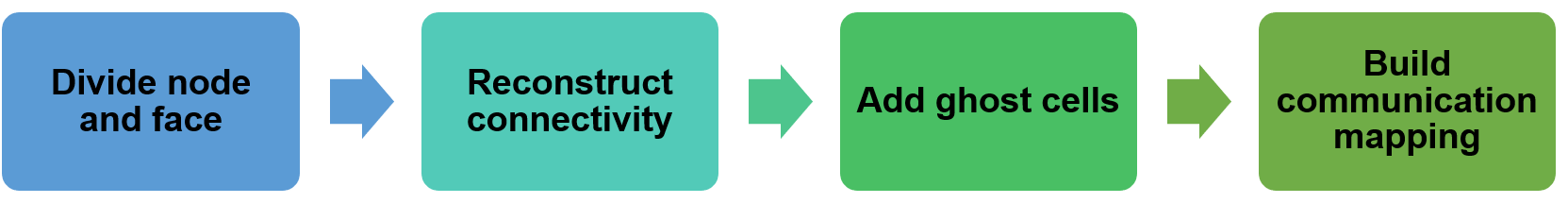}
\caption{\label{unm-geo-flowchart} Pre-processing of generating geometric information for each subdomain.}
\end{figure}

\subsection{Domain decomposition}
METIS library \cite{METIS-Kar} is used for the domain decomposition
of unstructured meshes in a reasonable way. For the load balance of
computation, the number of meshes assigned to each device is nearly
identical. For the efficiency of data communication, the number of
adjacent subdomain belong to different devices is minimized. METIS
only provides  the indices of cells for each subdomain, and the
geometric connectivity are broken after domain decomposition. The
pre-processes for data localization and independence are shown in
Figure.\ref{unm-geo-flowchart}.

\begin{figure}[!h]
\centering
\includegraphics[width=0.45\textwidth]{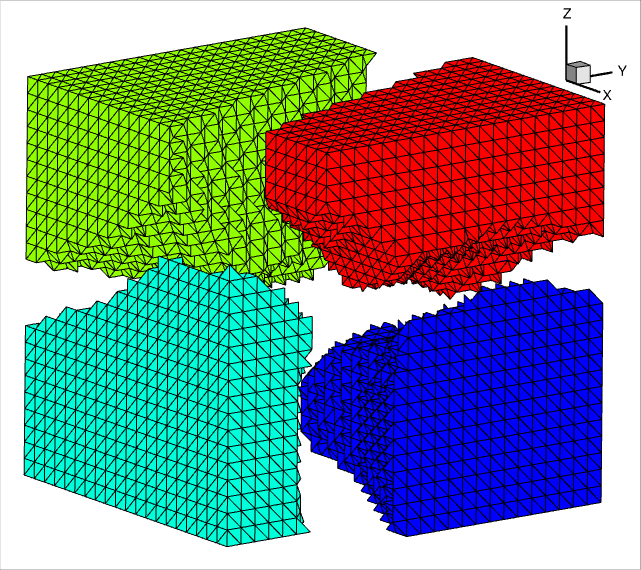}
\includegraphics[width=0.45\textwidth]{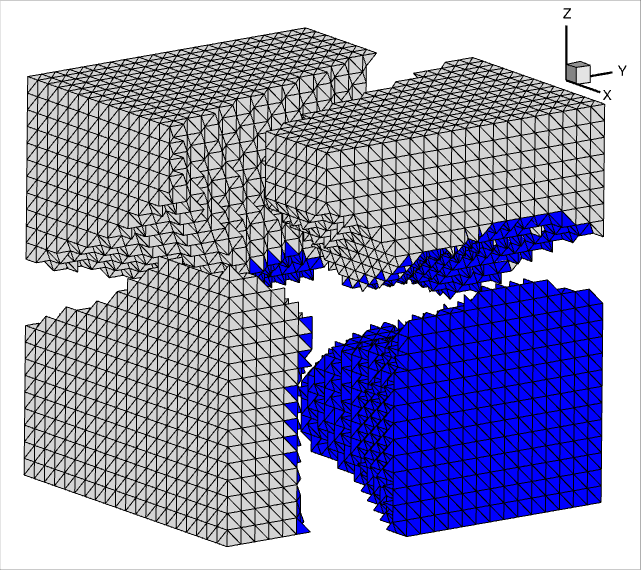}
\includegraphics[width=0.45\textwidth]{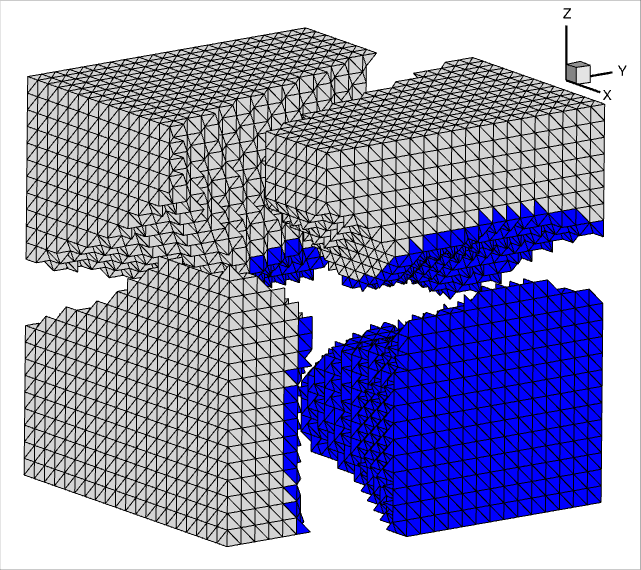}
\includegraphics[width=0.45\textwidth]{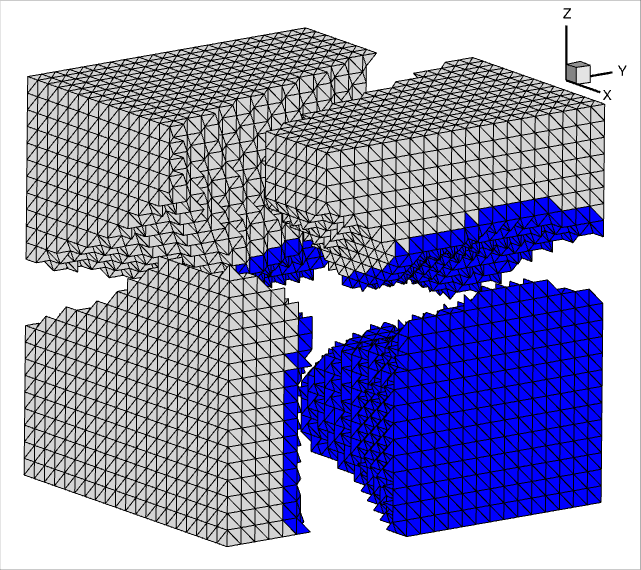}
\caption{\label{fig:unm-layer} Unstructured meshes are split by METIS and three-layer neighbors of subdomain in Blue distributed in other subdomains. }
\end{figure}

According to the cell information provided by METIS and the
geometric connectivities of the whole domain, the nodes and faces
can be classified into the corresponding subdomain. Duplicate nodes
and faces are included in different subdomains as well. A new set of
indices for nodes, faces, and cells is formed for each subdomain
after decomposition. There are two indices for nodes, faces and
cells, i.e.,  the old rank in the domain and the new rank in the
subdomains. The mappings of these indices are built by
two-dimensional arrays, i.e. \textit{SNodeToTNode} and
\textit{TNodeToSNode}. \textit{SNodeToTNode}[$i$][$j$] provides the
index of node $i$ of subdomain $j$ in the domain, and
\textit{TNodeToSNode}[$i$][$j$] offers the index of node $i$ of the
domain in subdomain $j$. The mappings of the indices of faces and
cells are established similarly. The connectivities of the domain
can be converted to the connectivities of each subdomain by using
these mappings, and the cell-to-face, face-to-node and cell-to-node
connectivities for each subdomain are also obtained. Due to large
stencil for WENO reconstruction, the data exchange becomes more
complicated. As shown in Figure.\ref{fig:unm-layer},  some
neighboring cells of each subdomain are located in other subdomains.
To compute the fluxes of inner boundary faces,  the geometric
information of the two cells sharing the inner boundary face and
their two-layer neighbors must be provided, which means three-layer
neighbors for the inner boundary should be accessed. These
three-layer neighbors are added to subdomains as ghost cells in the
pre-processing. According to the cell-to-neighbor connectivity,
these three-layer neighbors can be found, and their indices will be
added to the mappings between the domain and subdomains. The indices
of their nodes and faces will be added to subdomains in the same
way. Except for geometric information, the conservative variables
stored in three-layer neighbors should be transferred to ghost cells
after the update of flow variables. The indices of these neighbors
are used to build mappings between two subdomains for MPI
communication. These mappings are stored in two-dimensional array,
i.e., \textit{CellSend}. \textit{CellSend}[$i$][$2$] provides the
index of cell $i$ in the subdomain and the target subdomain. Through
such mappings, the conservative variables on the ghost cells can be
updated accurately. With the pre-processing above, the data
localization and independence of each subdomain are maintained.

\begin{figure}[!h]
\centering
\includegraphics[width=0.8\textwidth]{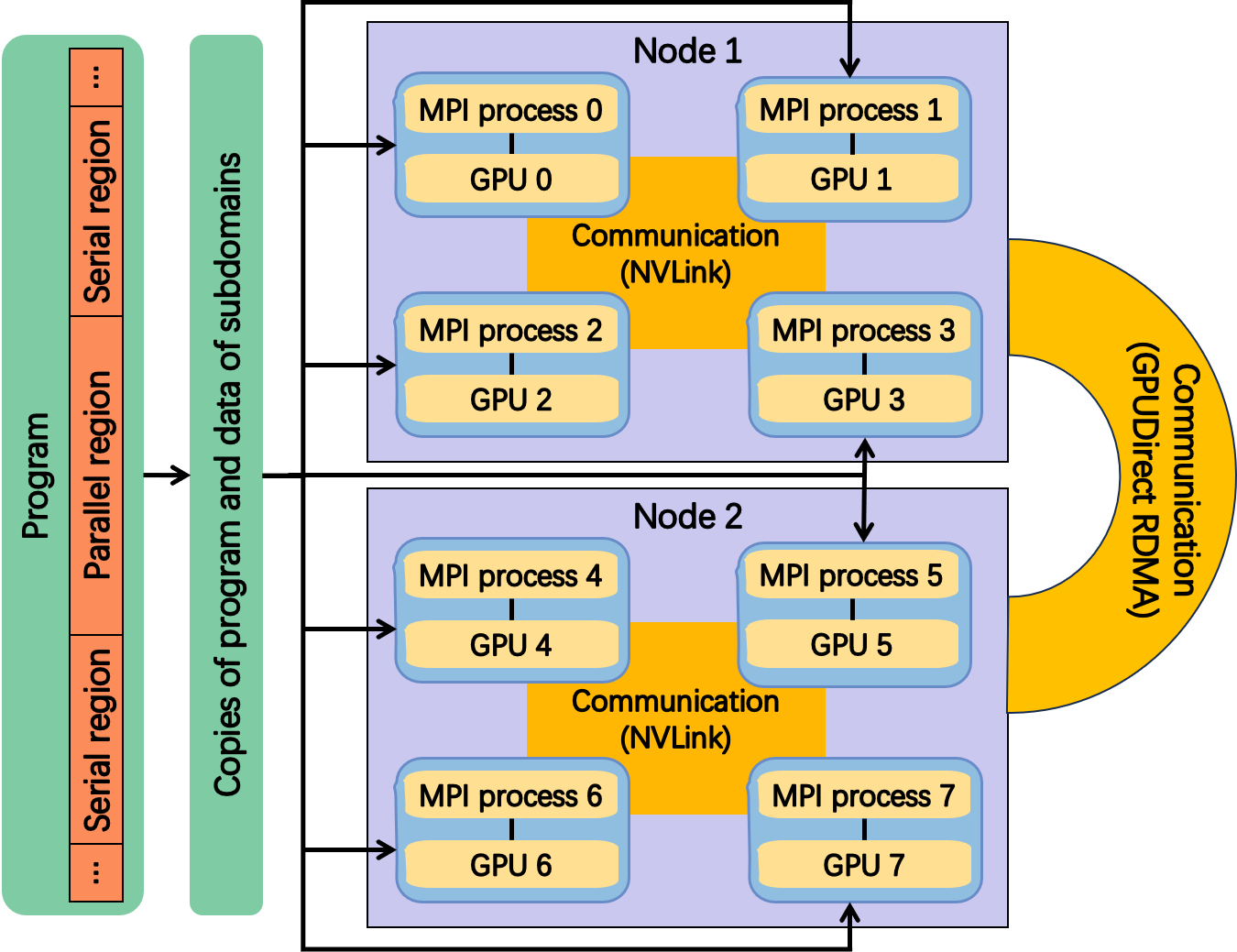}
\caption{\label{unm-mpi-gpu} Framework of multi-GPU with CUDA and MPI.}
\end{figure}

\subsection{Combination of MPI and CUDA}
MPI is a standard parallelization library adopted on distributed
memory platforms. It is widely used in CFD to meet the demands of
large-scale computations. The idea of creating MPI processes for
one-to-one GPU management and communication is maintained. The
framework of multiple-GPU implementation with CUDA and MPI is shown
in Figure.\ref{unm-mpi-gpu}. MPI is based on Single Program Multiple
Data (SPMD) model, in which a single program is executed
simultaneously with multiple data. For HGKS code with CUDA and MPI,
the copy of program and data of each subdomain are distributed to
each MPI process. The CUDA code for single GPU can be inherited as a
main part of multiple-GPU implementation, and the rest part for
multiple-GPU implementation is MPI communication. Typically, a
multi-GPU cluster consists of multiple nodes and each node contains
the same number of GPUs and CPUs. The CUDA-Aware MPI
library\cite{CUDA-MPI} is used in this study. It is an MPI
implementation that enables direct communication with GPU memory.
This capability markedly enhances parallel computing efficiency in
multi-GPU environments. NVLink and GPU Direct Remote Direct Memory
Access (RDMA)  are used to support direct GPU memory transfer within
nodes and across nodes, which improves communication efficiency
among GPUs.

\begin{figure}[!h]
\centering
\includegraphics[width=0.8\textwidth]{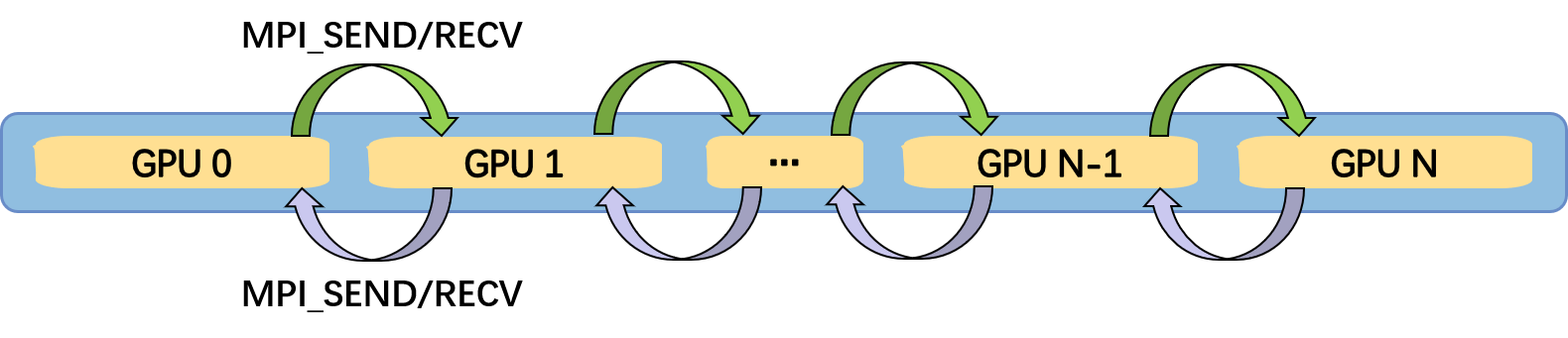}
\caption{\label{unm-s-mpi} Staggered blocking communication for case of structured meshes.}
\end{figure}

\begin{figure}[!h]
\centering
\includegraphics[width=0.7\textwidth]{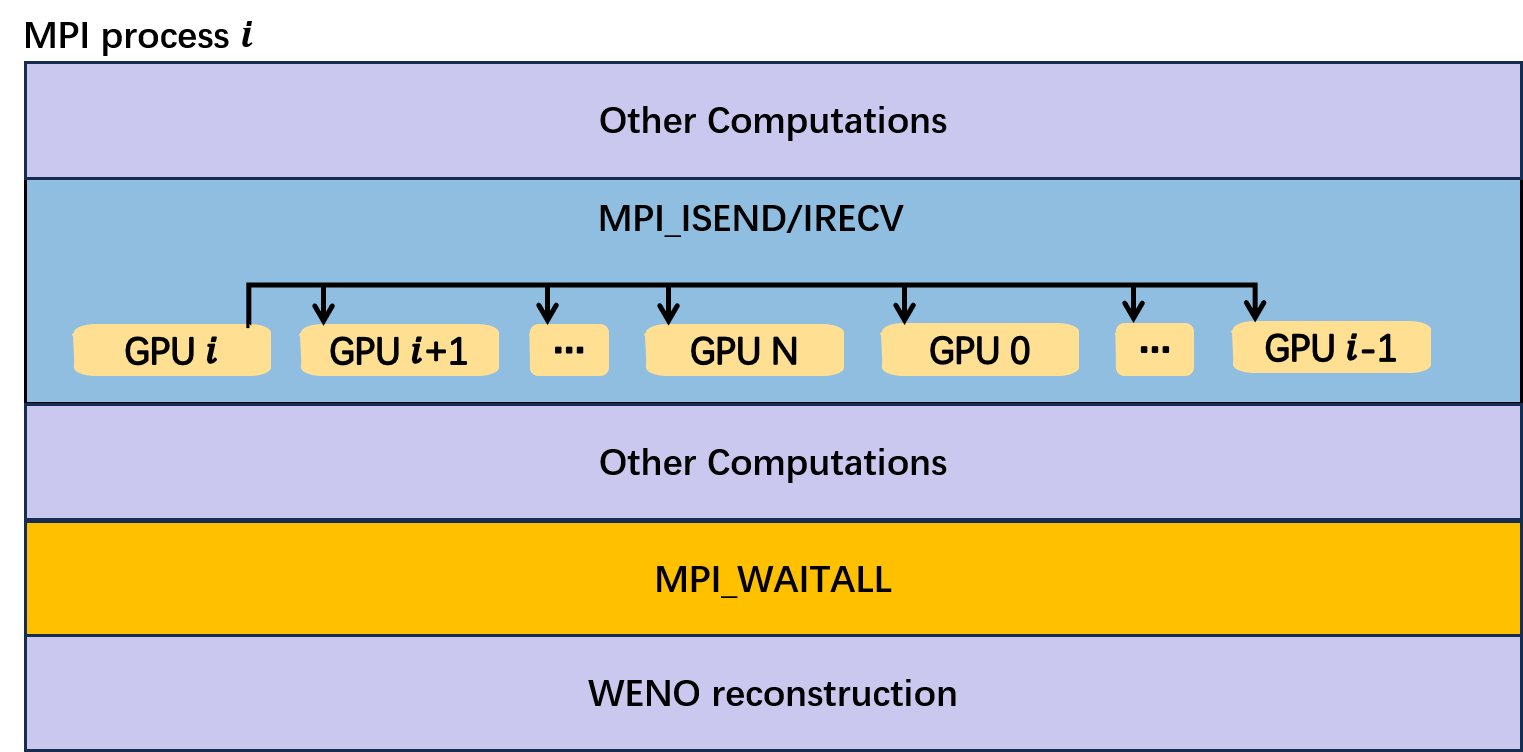}
\caption{\label{unm-mpi} Non-blocking communication for case of structured meshes.}
\end{figure}

In the previous study, as shown in Figure.\ref{unm-s-mpi}, the
blocking communication is employed in a staggered manner for the
HGKS code on structured meshes \cite{GKS-GPU-2}. Each subdomain only
need to transfer data with two neighboring subdomains, and the
communications in purple has to wait for the communications in
green. All other computations have to wait for the communication,
even though some computations do not require the data from the
communication. For the HGKS code on unstructured meshes, the
conservative variables of neighboring cells located in other
subdomains should be transferred to the target subdomain for WENO
reconstruction. Each subdomain has to communicate with all other
subdomains for the unstructured meshes. If the strategy for
structured mesh is extended to unstructured mesh,  a significant
amount of time will be spent for communication, and a new strategy
for communication must be proposed.  Non-blocking communication is
chosen for HGKS on unstructured meshes, and the brief setting of
communication is shown in Figure.\ref{unm-mpi}. Non-blocking
communication MPI$\_$ISEND/IRECV allows processes to initiate
communication operations and then proceed with other computations.
GPU $i$ can transfer data of three-layer neighbors with other GPUs
by using MPI$\_$ISEND/IRECV, simultaneously. The computations that
do not require these data can be executed on GPU during
communication. Before WENO reconstruction, MPI$\_$WAITALL is called
to ensure that all communication about GPU $i$ has been completed.
MPI$\_$ALLREDUCE and MPI$\_$BARRIER  are used to deal with the
operation of global unification. The indices of cells are
discontinuous, and the data packaging is added as a new
computational task in the program. By using the frame described
above, HGKS on unstructured meshes can be computed in parallel on
GPUs.

\section{Numerical tests}
In this section, the numerical examples will be presented to
validate the performance of GPU code by the accuracy test and the
flows passing through a sphere. For the GPU computations, the Nvidia
RTX A5000 GPU and Nvidia Tesla V100 GPU are used with Nvidia CUDA
and CUDA-Aware MPI of NVIDIA HPC SDK 21.7. As reference, the CPU
codes run on the station with 2 Intel(R) Xeon(R) Gold 5120 CPU using
OpenMP directives. The detailed parameters of CPU and GPU are given
in Table.\ref{GPU-CPU-A} and Table.\ref{GPU-CPU-B}, respectively.
For RTX A5000  GPU,  the GPU-GPU communication is achieved by
connection traversing PCIe, and there are $2$ RTX A5000 GPUs in one
node. For Tesla V100 GPU, there are $8$ GPUs inside one GPU node,
and more nodes are needed for more than $8$ GPUs. The GPU-GPU
communication in one GPU node is achieved by Nvidia NVLink. The
communication across GPU nodes can be achieved by GPU Direct RDMA
via iWARP, RDMA over Converged Ethernet (RoCE) or InfiniBand. In
this paper, RoCE is used for communication across GPU nodes.

\begin{table}[!h]
\begin{center}
\def\temptablewidth{0.95\textwidth}
 {\rule{\temptablewidth}{1.0pt}}
        \begin{tabular*}{\temptablewidth}{@{\extracolsep{\fill}}c|c|c|c}
            ~       & CPU                            & Clock rate & Memory size         \\
            \hline
            Station (56 cores) & Intel(R) Xeon(R) Gold 5120 CPU & $2.20$ GHz & $768$ GB \\
        \end{tabular*}
        {\rule{\temptablewidth}{1.0pt}}
        \caption{\label{GPU-CPU-A} The detailed parameters of CPU.}
    \end{center}
    \begin{center}
        \def\temptablewidth{0.95\textwidth}
        {\rule{\temptablewidth}{1.0pt}}
        \begin{tabular*}{\temptablewidth}{@{\extracolsep{\fill}}c|c|c}
            ~                            & Nvidia RTX A5000 & Nvidia Tesla V100 \\
            \hline
            Clock rate                   & 1.17 GHz         & 1.53 GHz          \\
            \hline
            Stream multiprocessor        & 64               & 80                \\
            \hline
            FP$32$ precision performance & 27.77 Tflops     & 15.7 Tflops       \\
            \hline
            FP$64$ precision performance & 867.8 Gflops     & 7834 Gflops       \\
            \hline
            GPU memory size              & 24 GB            & 32 GB             \\
            \hline
            Memory bandwidth             & 768 GB/s         & 897 GB/s          \\
        \end{tabular*}
        {\rule{\temptablewidth}{1.0pt}}
        \caption{\label{GPU-CPU-B} The detailed parameters of GPU.}
    \end{center}
\end{table}

\begin{figure}[!h]
    \centering
    \includegraphics[width=0.5\textwidth]{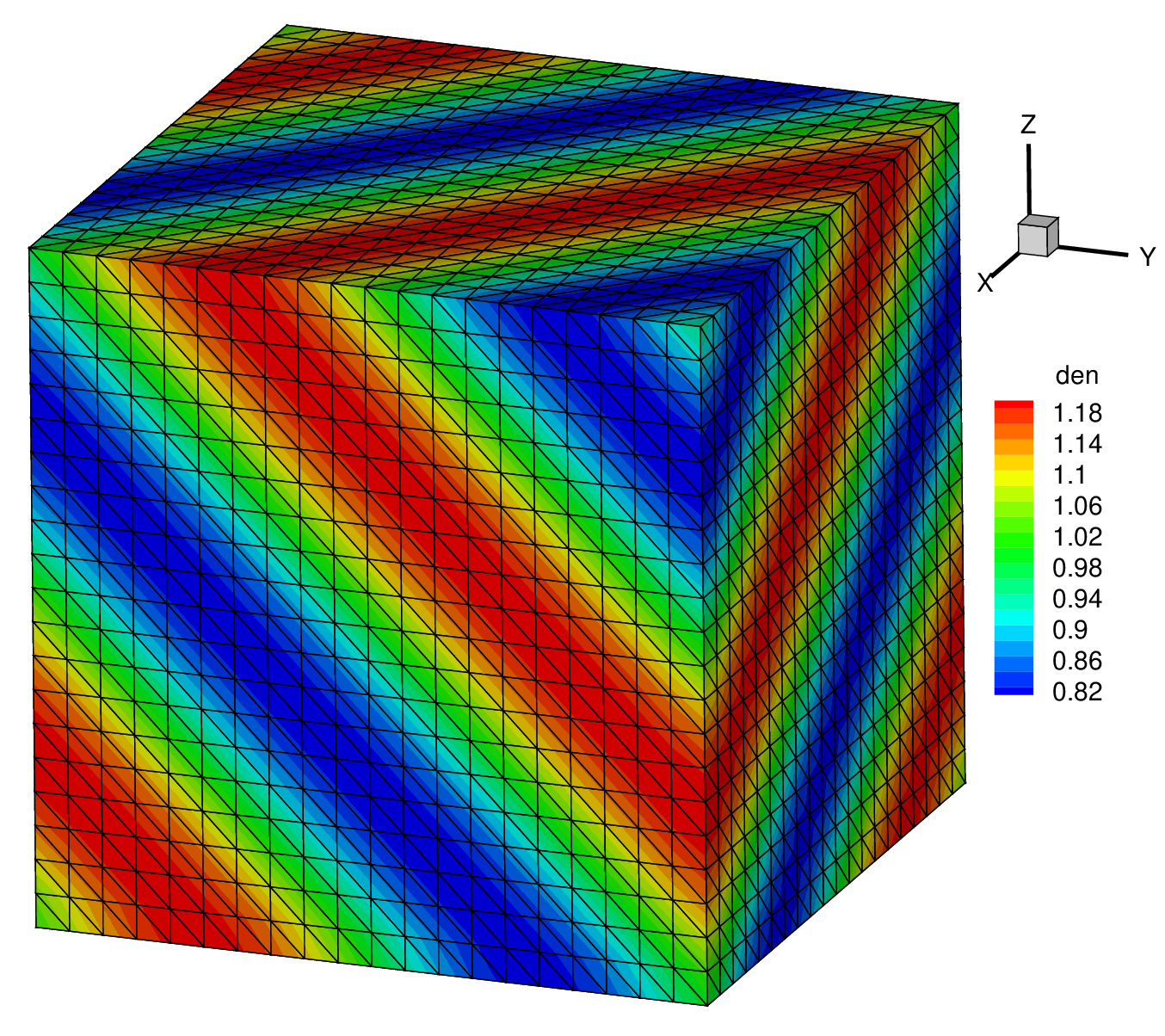}
    \caption{\label{accuracy-mesh} Accuracy test: the computational mesh with $20^3\times 6$ cells and the density distribution at $t=2$.}
\end{figure}

\subsection{Accuracy test}
In this case, the three-dimensional advection of  density
perturbation is used to test the order of accuracy with the
tetrahedron meshes. For this case, the computational domain is
$[0,2]\times[0,2]\times[0,2]$ and the initial condition is given as
follows
\begin{align*}
\rho_0 & (x, y, z)=1+0.2\sin(\pi(x+y+z)),~p_0(x,y,z)=1, \\
           & U_0(x,y,z)=1,~V_0(x,y,z)=1,~W_0(x,y,z)=1.
\end{align*}
The periodic boundary condition is applied on all boundaries, and the
exact solution is
\begin{align*}
\rho(x,y & ,z,t)=1+0.2\sin(\pi(x+y+z-t)),~p(x,y,z,t)=1, \\
& U(x,y,z,t)=1,~V(x,y,z,t)=1,~W(x,y,z,t)=1.
\end{align*}
In smooth flow regions, $\tau=0$ and the gas-distribution function
reduces to
\begin{align*}
    f(\boldsymbol{x}_{G},t,\boldsymbol{u},\xi)=g(1+At).
\end{align*}
For this case, a series of meshes with $6 \times N^3$ tetrahedron
cells are used, where every cubic is divided into six tetrahedron
cells.  The computational mesh with $20^3\times 6$ cells are shown
in Figure. \ref{accuracy-mesh}. The density distribution at $t=2$ is
also shown in Figure. \ref{accuracy-mesh}, and the $L^1$ and $L^2$
errors and orders of accuracy at $t=2$ are presented in
Tab.\ref{GPU-CPU-TET-ACC}, where the expected order of accuracy is
achieved.

\begin{table}[!h]
    \begin{center}
        \def\temptablewidth{0.85\textwidth}
        {\rule{\temptablewidth}{1.0pt}}
        \begin{tabular*}{\temptablewidth}{@{\extracolsep{\fill}}c|c|c|c|c}
            mesh           & $L^1$ error & Order  & $L^2$ error & Order  \\
            \hline
            $10^3\times 6$ & 6.6070E-02  & ~~     & 2.5938E-02  & ~~     \\
            \hline
            $20^3\times 6$ & 8.7117E-03  & 2.9230 & 3.4159E-03  & 2.9247 \\
            \hline
            $40^3\times 6$ & 1.0994E-03  & 2.9862 & 4.3094E-04  & 2.9867 \\
            \hline
            $80^3\times 6$ & 1.3768E-04  & 2.9973 & 5.3968E-05  & 2.9973 \\
            \hline
            $160^3\times6$ & 1.7252E-05  & 2.9965 & 6.7761E-06  & 2.9936 \\
        \end{tabular*}
        {\rule{\temptablewidth}{1.0pt}}
        \caption{\label{GPU-CPU-TET-ACC} Accuracy tests: errors and orders of
            accuracy with tetrahedron meshes.}
    \end{center}
\end{table}

\begin{table}[!h]
    \begin{center}
        \def\temptablewidth{0.85\textwidth}
        {\rule{\temptablewidth}{1.0pt}}
        \begin{tabular*}{\temptablewidth}{@{\extracolsep{\fill}}c|c|c|c|c|c}
            Mesh size      & Station & RTX A5000 & Speedup & Tesla V100 & Speedup \\
            \hline
            $6\times 10^3$ & 5.44       & 1.05         & 5.18     & 0.55             & 9.89    \\
            \hline
            $6\times 20^3$ & 58.12      & 10.1         & 5.75     & 5.99             & 9.70    \\
            \hline
            $6\times 40^3$ & 840.89     & 150.57       & 5.58     & 86.01            & 9.78    \\
            \hline
            $6\times 80^3$ & 12821.05   & 2342.37      & 5.47     & 1333.33          & 9.62     \\
        \end{tabular*}
        {\rule{\temptablewidth}{1.0pt}}
        \caption{\label{GPU-CPU-TG-C} Accuracy tests: The detailed computational times(s) and speedup for
RTX A5000, and  56-core station CPU with OpenMP is used as
reference.}
    \end{center}
    \begin{center}
        \def\temptablewidth{0.55\textwidth}
        {\rule{\temptablewidth}{1.0pt}}
        \begin{tabular*}{\temptablewidth}{@{\extracolsep{\fill}}c|c|c}
            No. GPUs  & $6\times40^3$ & Parallel efficiency  \\
            \hline
            1                             & 86.01             & ~         \\
            \hline
            2                             & 44.26             & 0.97      \\
            \hline
            4                             & 23.79             & 0.90       \\
            \hline
            8                             & 13.73             & 0.78     \\
            \hline
            \hline
            No. GPUs  & $6\times80^3$ & Parallel efficiency   \\
            \hline
            1                                & 1333.33           & ~      \\
            \hline
            2                              & 672.23            & 0.99     \\
            \hline
            4                             & 343.87            & 0.97       \\
            \hline
            8                        & 180.19            & 0.92      \\
            \hline
            \hline
            No. GPUs  & $6\times160^3$ & Parallel efficiency \\
            \hline
            4                   & 5375.23           & ~       \\
            \hline
            8             & 2710.12           & 0.99
        \end{tabular*}
        {\rule{\temptablewidth}{1.0pt}}
        \caption{\label{GPU-CPU-TG-D}  Accuracy tests:  The detailed computational time(s) and parallel efficiency for Tesla V100.}
    \end{center}
\end{table}

To test the performance of single-GPU code, this case is run with
both Nvidia RTX A5000 and Nvidia Tesla V100 GPUs using Nvidia CUDA.
As a comparison , the CPU code running on the Intel Xeon(R) Gold
5120 CPU is also tested. The execution times with different meshes
are shown in Table.\ref{GPU-CPU-TG-C}, where the total execution
time for CPU and GPUs are given at $t=2$. Due to the limitation of
single GPU memory size, the tetrahedron meshes from $6\times10^3$ to
$6\times80^3$ cells are used. For the CPU computation, $56$ cores
are utilized with OpenMP parallel computation, and the corresponding
execution time is used as the base for following comparisons. The
speedups  are  given in Table.\ref{GPU-CPU-TG-C}, which is defined
as
\begin{align*}
\displaystyle\text{Speedup}=\frac{\text{CPU time}}{\text{GPU time}}.
\end{align*}
Compared with the CPU code using 56 cores, $\approx 5$x speedup is
achieved by single RTX A5000 GPU, and $\approx 9$x speedup is
achieved by single Tesla V100 GPU. Even though the Tesla V100 is
approximately $15$ times faster than the RTX A5000 in FP64 precision
computation ability, Tesla V100 only performs $2$ times faster than
RTX A5000. Due to the limitation of single-GPU, the multiple-GPU
accelerated HGKS code has to be implemented. To further show the
performance of GPU computation, the scalability is defined as
\begin{equation*}
\displaystyle\text{Parallel efficiency}=\frac{\text{single-GPU time}}{\text{multiple-GPU time}\times \text{number of GPUs}}.
\end{equation*}
The detailed data for Tesla V100 GPU are given in
Table.\ref{GPU-CPU-TG-D}, in which $1$ to $8$ GPUs are used for the
case with $6\times40^3$ and $6\times80^3$ cells, and $4$ to $8$ GPUs
are used for the case with $6\times160^3$ cells. For the ideal
parallel computations, the parallel efficiency would be equal to $1$
with the increase of GPUs. However, such parallel efficiency is not
possible due to the communication delay among the computational
cores and the idle time of computational nodes associated with load
balancing. As shown in Table.\ref{GPU-CPU-TG-D}, the explicit
formulation of HGKS scales properly with the increasing number of
GPU. Each GPU in multiple-GPU implementation can achieve more than
97$\%$ of the efficiency of single-GPU implementation with suitable
computational work-load.

\begin{table}[!h]
    \begin{center}
        \def\temptablewidth{0.95\textwidth}
        {\rule{\temptablewidth}{1.0pt}}
        \begin{tabular*}{\temptablewidth}{@{\extracolsep{\fill}}c|c|c|c}
        RTX A5000 & Computation time with FP32 & Computation time with FP64 & Ratio \\
        \hline
        $6\times40^3$      & 32.53                         &    150.57                      & 4.63    \\
        \hline
        $6\times80^3$      & 521.91                         & 2342.37                         & 4.49    \\
        \hline
        RTX A5000 & Memory cost with FP32         & Memory cost with FP64         & Ratio \\
        \hline
        $6\times40^3$     & 1470                         & 2578                         & 1.75    \\
        \hline
        $6\times80^3$      & 7762                        & 14822                        & 1.91    \\
        \hline
        \hline
        Tesla V100 & Computation time with FP32 & Computation time with FP64 & Ratio \\
        \hline
        $6\times40^3$       & 32.78                         & 86.01                         & 2.62    \\
        \hline
        $6\times80^3$      & 490.64                         & 1333.33                         & 2.72    \\
        \hline
        Tesla V100 & Memory cost with FP32         & Memory cost with FP64         & Ratio \\
        \hline
        $6\times40^3$     & 1843                        & 3157                         & 1.71    \\
        \hline
        $6\times80^3$      & 8095                         & 14993                       & 1.85    \\
        \end{tabular*}
        {\rule{\temptablewidth}{1.0pt}}
    \caption{\label{accuracy-precision} Accuracy test: the
        comparison of computational time (s) and memory (MiB) for FP32 and FP64
        precision with single RTX A5000 and Tesla V100.}
    \end{center}
    \begin{center}
        \def\temptablewidth{0.85\textwidth}
        {\rule{\temptablewidth}{1.0pt}}
        \begin{tabular*}{\temptablewidth}{@{\extracolsep{\fill}}c|cc|cc}
            mesh     & $L^1$ error  &    Order      &  $L^2$ error &  Order   \\
            \hline
            $6\times10^3$ & 6.6074E-02  & ~      & 2.5940E-02  & ~        \\
            \hline
            $6\times20^3$ & 8.7154E-03  & 2.9225 & 3.4173E-03  & 2.9242 \\
            \hline
            $6\times40^3$ & 1.1019E-03  & 2.9835 & 4.3188E-04  & 2.9842 \\
            \hline
            $6\times80^3$ & 1.4119E-04  & 2.9644 & 5.5303E-05  & 2.9652 \\
        \end{tabular*}
        {\rule{\temptablewidth}{1.0pt}}
    \end{center}
    \caption{\label{tab-3d-Accuracy-fp32} Accuracy test: errors and orders of accuracy with tetrahedron cells using RTX A5000 and FP32 precision.}
\end{table}

In this case, the GPU-accelerated HGKS is compiled with both FP32
precision and FP64 precision. For most GPUs, the  performance of
FP32  precision is stronger than FP64 precision. Because of the
reduction in device memory and improvement of arithmetic
capabilities on GPUs, the benefits can be achieved by using FP32
precision compared to FP64 precision. In view of these strength,
FP32-based and mixed-precision-based high-performance computing
start to be explored \cite{GPU-9,GPU-10}. The
comparison study with FP32 and FP64 precision for the direct
numerical simulation is provided in our previous work
\cite{GKS-GPU-2}. For the HGKS code on unstructured mesh, such
comparison is also given. For simplicity, the memory cost and
computational time for the accuracy test are provided in
Table.\ref{accuracy-precision}, where the execution times are given
in terms of second and the memory cost is in MiB.  As expected, the
memory of  FP32 precision is about half of that of FP64 precision
for both RTX A5000 and Tesla V100 GPUs. Compared with FP64-based
simulation, 4.6x speedup is achieved for RTX A5000 GPU and 2.6x
speedup is achieved for Tesla V100 GPU. The errors and orders of
accuracy with tetrahedron cells using RTX A5000 and FP32 precision
is given in Table.\ref{tab-3d-Accuracy-fp32}, which deviates with
the error in Table.\ref{GPU-CPU-TG-C} slightly.  In terms of the
problems without very strict requirements in accuracy, such as the
third-order scheme on unstructured meshes, FP32 precision may be
used due to its improvement of efficiency and reduction of memory.
However, for the direct numerical simulation of turbulences,  the
numerical results indicate that FP32 precision is not enough
\cite{GKS-GPU-2}. It strongly suggests that the FP64 precision
performance of GPU still requires to be improved to accommodate the
increasing requirements of GPU-based HPC.

\begin{figure}[!h]
\centering
\includegraphics[width=0.55\textwidth]{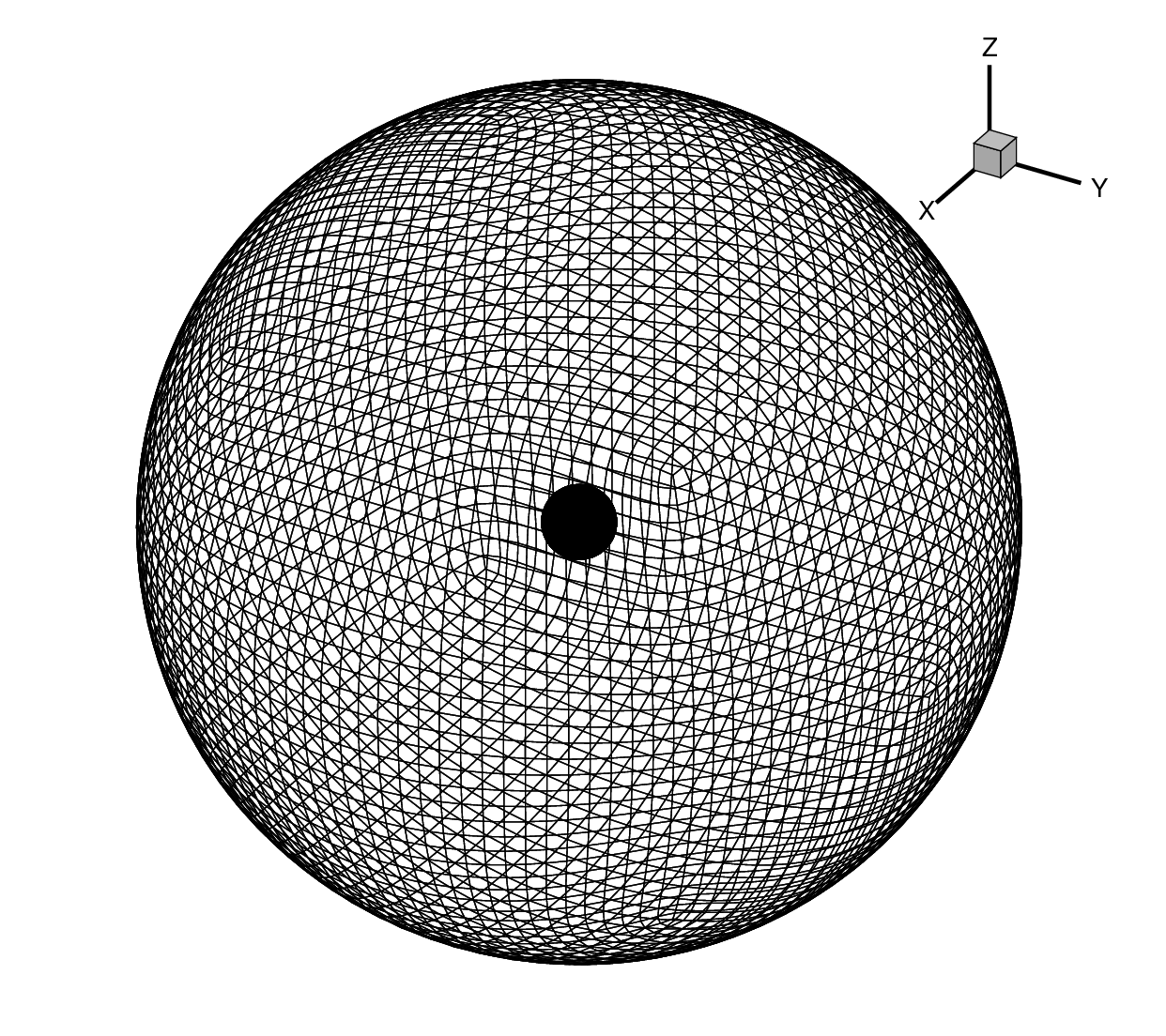}
\caption{\label{sphere-mesh} Flows passing through a sphere: the mesh distribution with $6\times16\times16\times32$ cells.}
\end{figure}

\subsection{Flows passing through a sphere}
\begin{figure}[!h]
\centering
\includegraphics[width=0.495\textwidth]{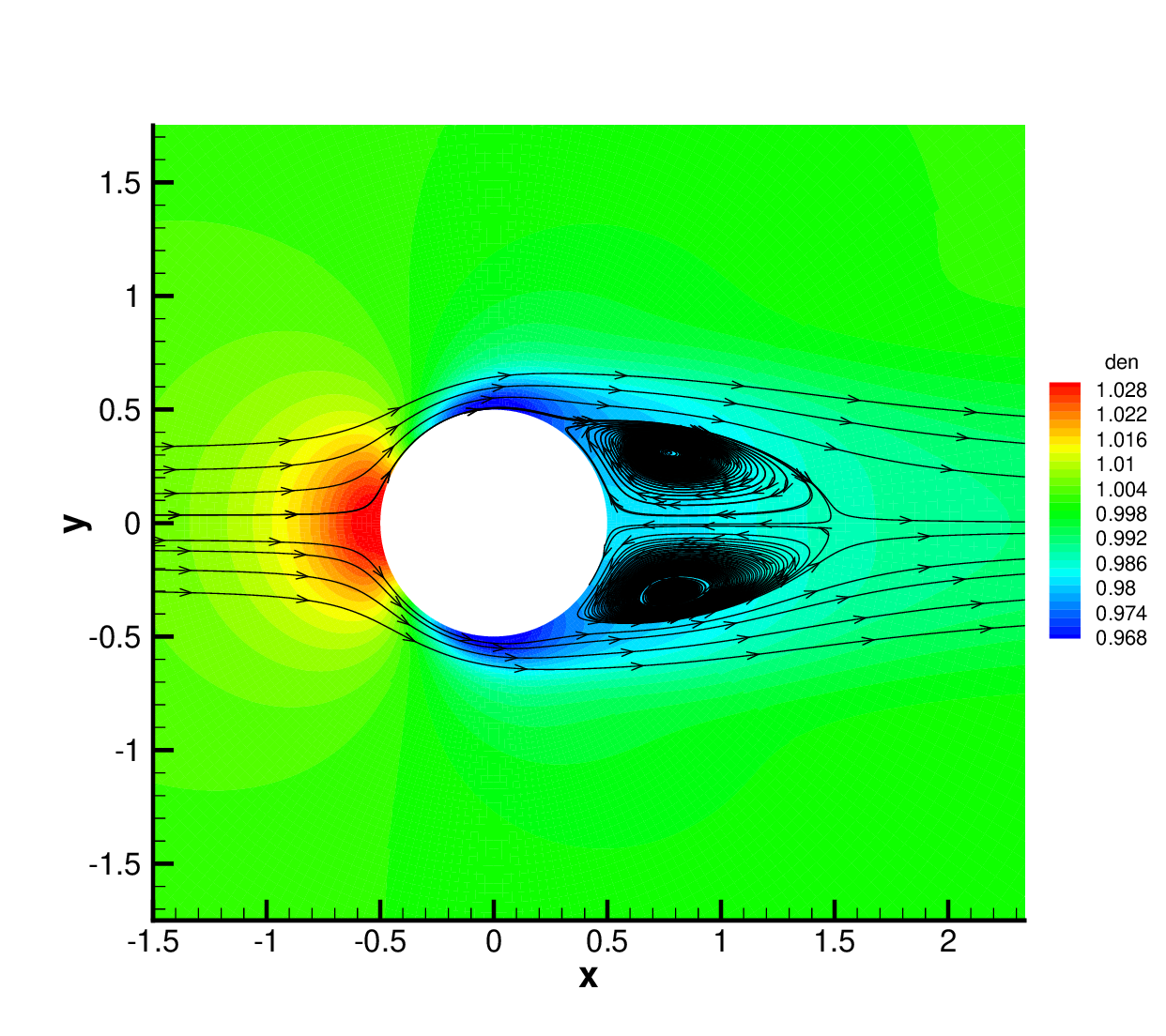}
\includegraphics[width=0.495\textwidth]{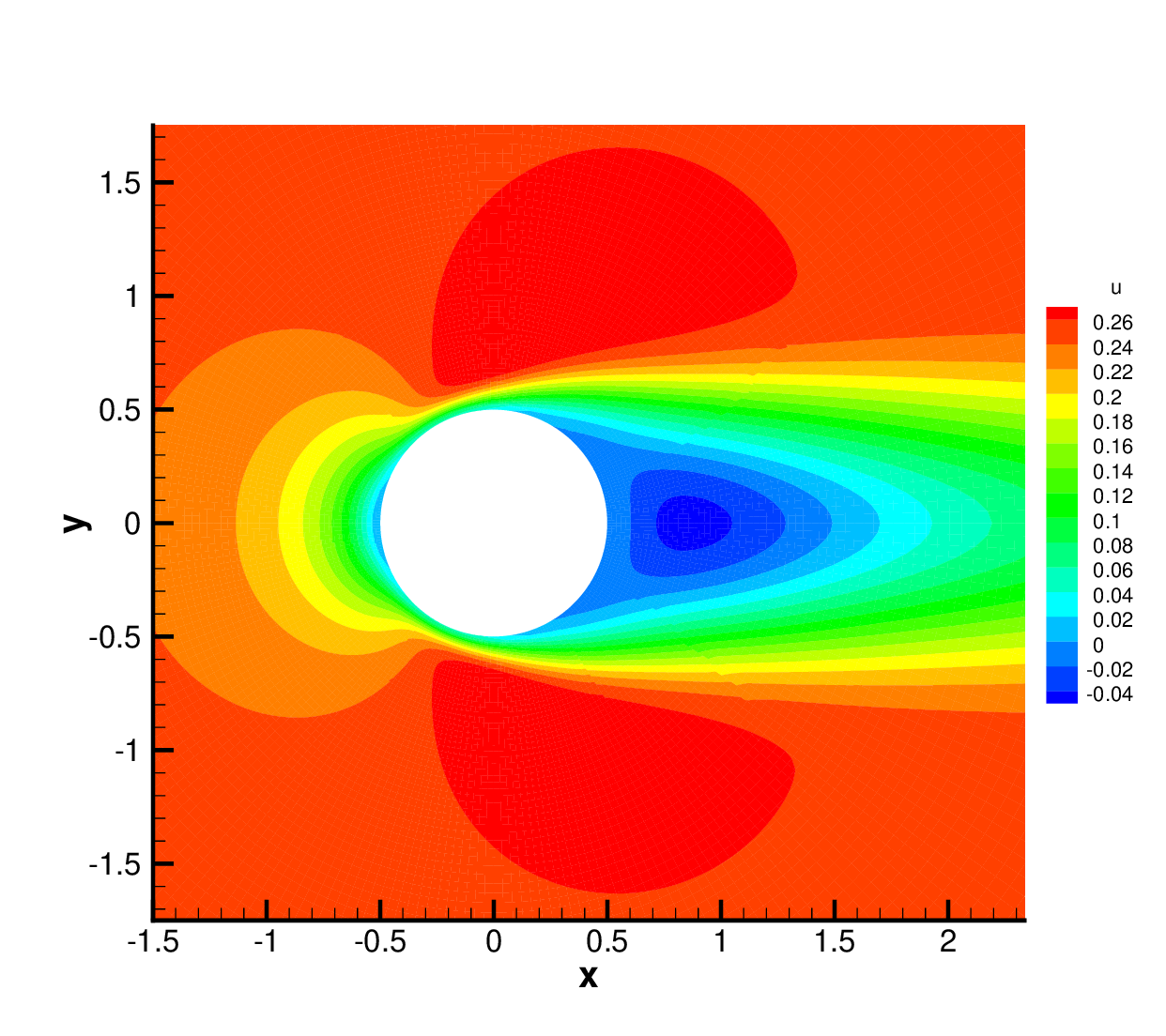}
\includegraphics[width=0.495\textwidth]{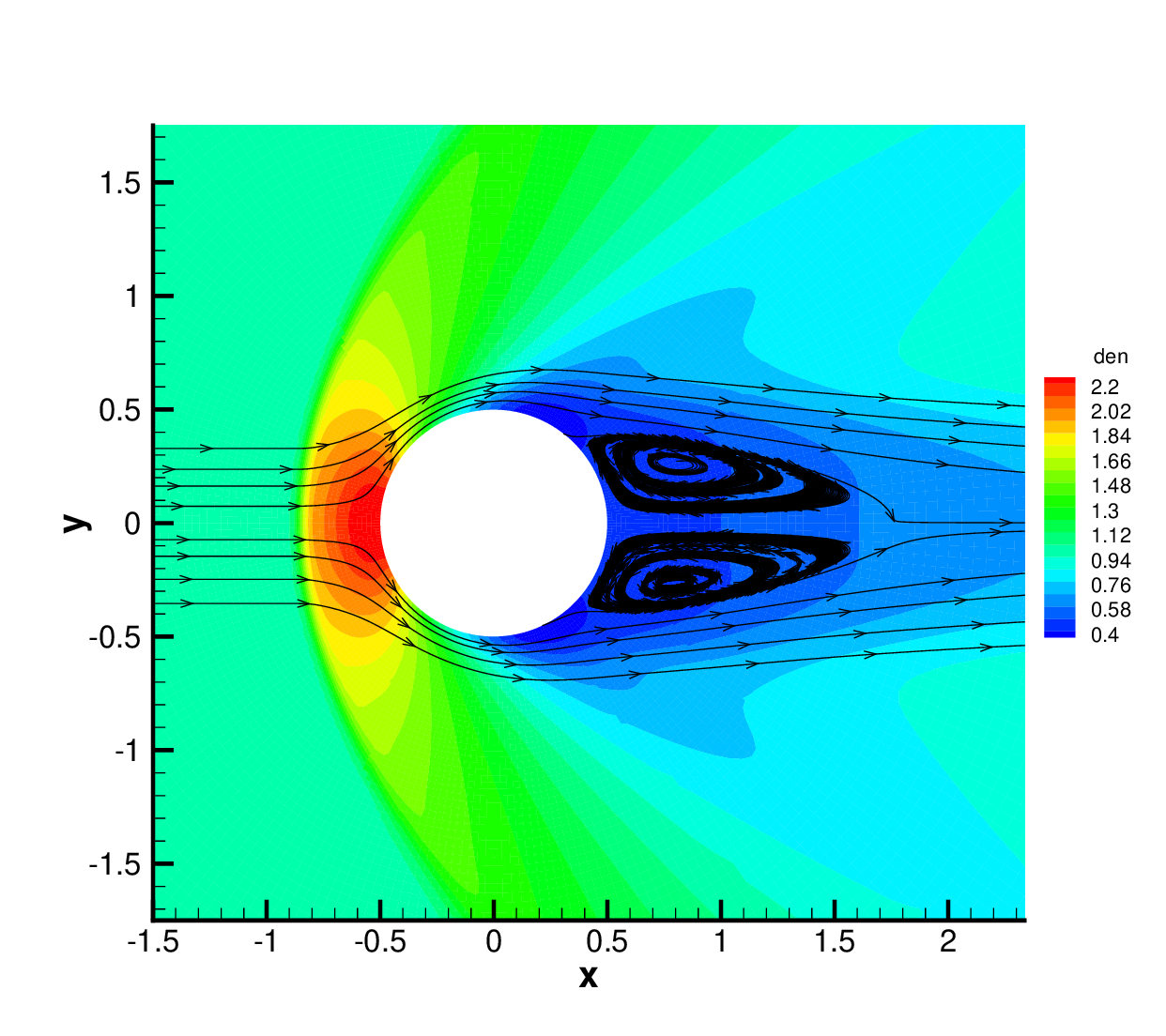}
\includegraphics[width=0.495\textwidth]{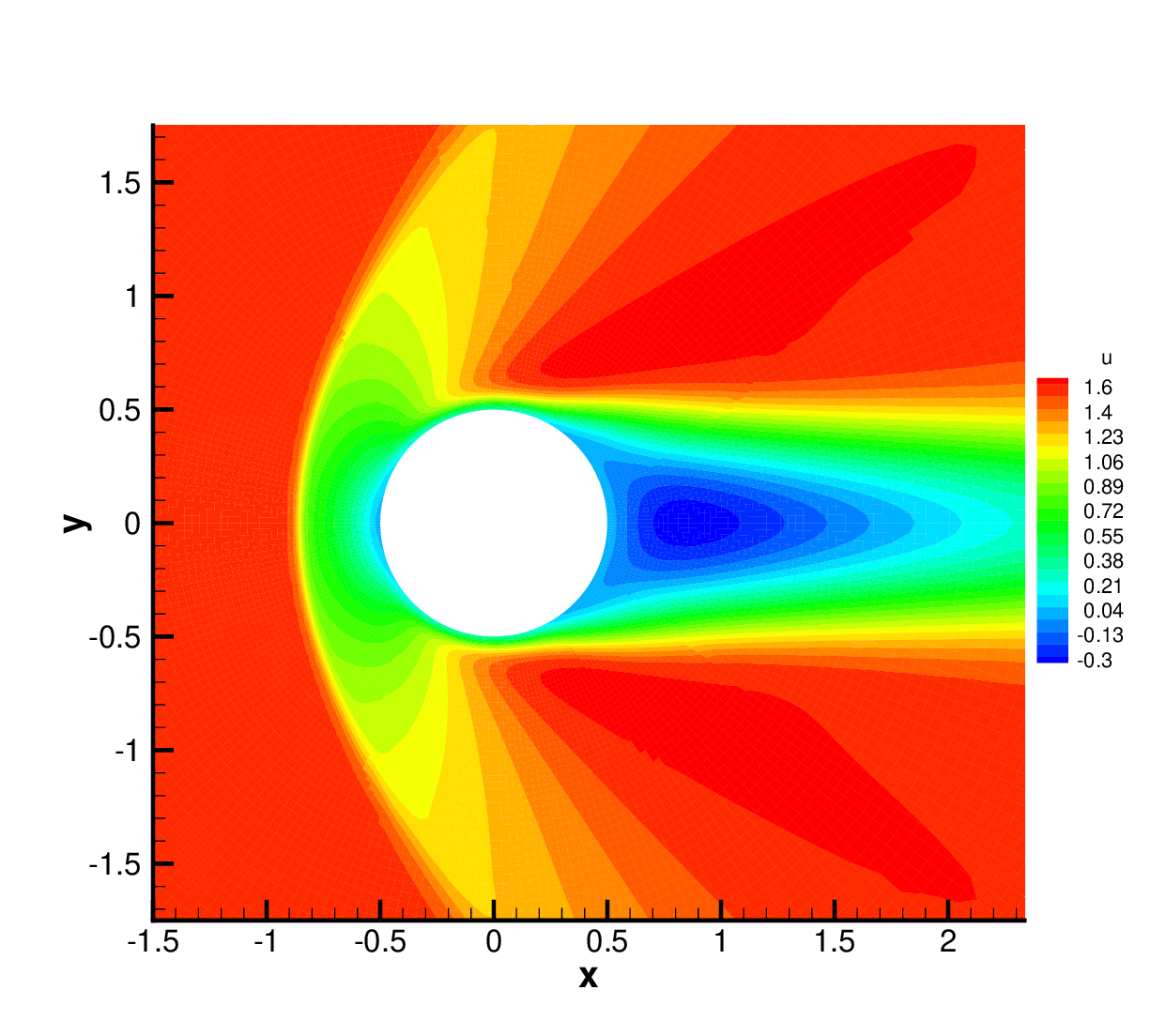}
\caption{\label{sphere-distribution} Flows passing through a sphere: the density, streamline and mach distribution for
the case with $Re=118$ and $Ma_\infty=0.2535$ (top) and  $Re=300$ and $Ma_\infty=1.5$ (bottom) using $12\times 64^3$ cells.}
\end{figure}

This case is used to test the capability in resolving from the
low-speed to transsonic flows, and the initial condition is given as
a free stream condition
\begin{align*}
(\rho,U,V,W,p)_{\infty} = (1, Ma_\infty,0,0, 1/\gamma),
\end{align*}
where $\gamma= 1.4$, $Ma_\infty$ is the Mach number of free stream, and the diameter of the sphere is $D=1$.
As shown in Figure.\ref{sphere-mesh}, this case is performed by the hexahedral meshes. These meshes are composed with
six parts and each part contains $N\times N\times 2N$ cells, and the total number of cell is $12\times N^3$.
The subsonic case with $Re=118$ and $Ma_\infty=0.2535$ and  the transsonic case with $Re=300$ and $Ma_\infty=1.5$ are  tested.
The inlet and outlet boundary conditions are given according to Riemann invariants, and the non-slip adiabatic boundary condition
is imposed for viscous flows on the surface of sphere.  The dynamic viscosity is given by
\begin{align*}
\mu=\mu_\infty(\frac{T}{T_\infty})^{0.7},
\end{align*}
where $T_\infty$ and $\mu_\infty$ are free stream temperature and
viscosity.  The density, velocity and streamline distributions at
vertical centerline planes are shown in
Figure.\ref{sphere-distribution} for the subsonic case and
transsonic case using $12\times 64^3$ cells. The quantitative
results of separation angle $\theta$  and closed wake length $L$ are
given in Table.\ref{sphere-2535-table} and
Table.\ref{sphere-090-table}, respectively. The current scheme gives
a better values of $L$ and $\theta$, which are closer to the
experiment data. The performance would be better with a refined
mesh. Quantitative results agree well with the benchmark solutions
\cite{Case-Nagata}, and the slight deviation of compact gas-kinetic
scheme \cite{GKS-high-2} might caused by the coarser mesh.

\begin{table}[!h]
\begin{center}
\def\temptablewidth{0.8\textwidth}{\rule{\temptablewidth}{0.9pt}}
\begin{tabular*}{\temptablewidth}{@{\extracolsep{\fill}}c|c|c|c}
Scheme &Computational Mesh &  $L$ & $\theta$ \\
\hline
Current result                                      & 3145728  cells, Hex        &   0.99  & 128.4\\
LUSGS GKS    \cite{GKS-high-4}                                   & 190464  cells, Hex          &   0.91  & 124.5\\
GMRES GKS   \cite{GKS-high-4}                                   & 190464  cells, Hex          &  0.91  & 124.5 \\
GMRES DDG  \cite{Case-Cheng}     & 1608680 cells, Hybrid     &  0.96   & 123.7 \\
Experiment \cite{Case-Taneda}                    & -                         & 1.07    & 151.0\\
\end{tabular*}
{\rule{\temptablewidth}{1.0pt}}
\end{center}
\caption{\label{sphere-2535-table} Flows  passing through a sphere:
the quantitative comparisons of closed wake length $L$ and
separation angle $\theta$ for $Re=118$ and $Ma_\infty=0.2535$.}
\begin{center}
\def\temptablewidth{0.8\textwidth}
{\rule{\temptablewidth}{1.0pt}}
\begin{tabular*}{\temptablewidth}{@{\extracolsep{\fill}}c|c|c|c}
scheme     & Computational Mesh                  &  $L$  &    $\theta$  \\
\hline
Current result                     &  3145728  cells, Hex   &  1.18 & 136.8    \\
LUSGS GKS   \cite{GKS-high-4}               &  190464  cells, Hex   &  1.17 & 135.1    \\
GMRES GKS   \cite{GKS-high-4}               &  190464  cells, Hex &  1.17 & 135.4
\end{tabular*}
{\rule{\temptablewidth}{1.0pt}}
 \end{center}
 \caption{\label{sphere-090-table} Flows  passing through a sphere: the quantitative comparisons of
closed wake length $L$ and separation angle $\theta$ for  $Re=300$
and $Ma_\infty=1.5$.}
\end{table}

\begin{table}[!h]
    \begin{center}
        \def\temptablewidth{0.75\textwidth}
        {\rule{\temptablewidth}{1.0pt}}
        \begin{tabular*}{\temptablewidth}{@{\extracolsep{\fill}}c|c|c}
            No. GPUs ($N=16$) &Total  &  Communication $(\%)$  \\
            \hline
            1              &  0.105               & ~              \\
            \hline
            2              & 0.056               & 0.426              \\
            \hline
            4              & 0.031               & 0.823                \\
            \hline
            8                & 0.019               & 1.648             \\
            \hline
            \hline
            No. GPUs  ($N=32$) &Total  & Communication $(\%)$    \\
            \hline
            1              &  1.747                  & ~            \\
            \hline
            2             & 0.883                 & 0.024             \\
            \hline
            4             &  0.452                 & 0.115           \\
            \hline
            8            & 0.232                 & 0.554            \\
            \hline
            \hline
            No. GPUs  ($N=64$) &Total   &  Communication $(\%)$     \\
            \hline
            1           &  29.096        & ~           \\
            \hline
            2          &14.587         & 0.034              \\
            \hline
            4         & 7.331          & 0.082              \\
            \hline
            8          &   3.688         & 0.149              \\
            \hline
            16                  & 1.929             & 3.169    \\
            \hline
            \hline
            No. GPUs   ($N=128$) &Total &  Communication $(\%)$    \\
            \hline
            8                &  60.058              & 0.063               \\
            \hline
            16                             & 30.876             & 2.183    \\
        \end{tabular*}
        {\rule{\temptablewidth}{1.0pt}}
        \caption{Flows passing through a sphere:  the detailed computational time (h) and parallel efficiency for Tesla V100.}
    \end{center}
\end{table}

\begin{figure}[!h]
    \centering
    \includegraphics[width=0.55\textwidth]{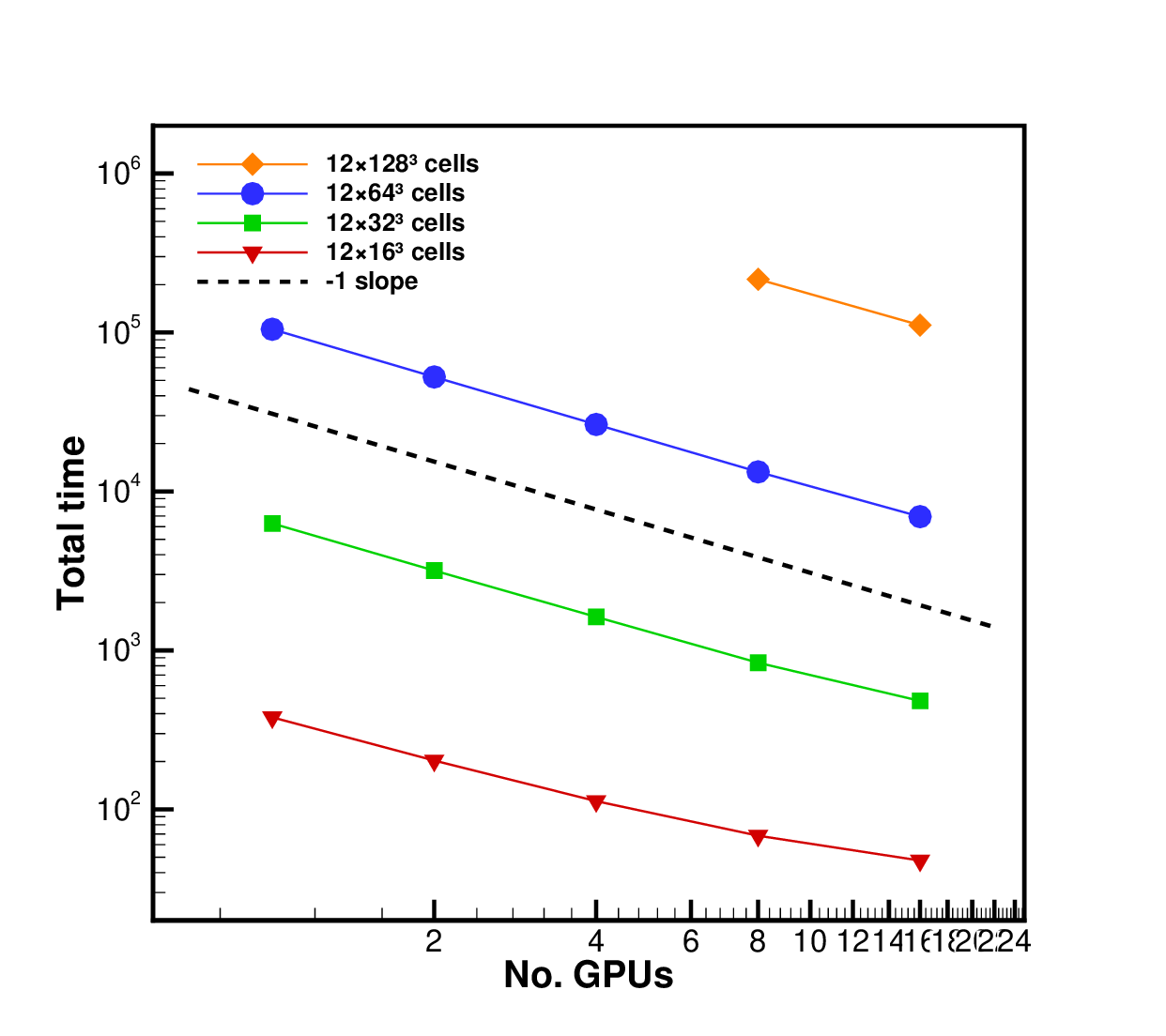}
    \caption{\label{mpi-efficience-hex} Flows passing through a sphere: the speedup of Tesla V100 GPUs with $12\times16^3$ to  $12\times128^3$ cells.}
\end{figure}

To show the performance with multiple GPUs, this case is tested by
the meshes from $12\times16^3$ to $12\times128^3$ cells. For the
mesh with $12\times16^3$ and $12\times32^3$ cells, the maximum of
GPU is 8, and NVLink  are used within node. For the mesh with
$12\times64^3$ and $12\times128^3$ cells, the maximum number  of GPU
is 16. NVLink and RDMA via RoCE are used within node and across
nodes respectively. The log-log plot for $n$  and total time are
shown in Figure.\ref{mpi-efficience-hex}. Conceptually, the total
computation amount increases with a factor of $16$, when the number
of $N$ doubles. With the log-log plot for $n$ and  total time, an
ideal scalability would follow $-1$ slope. However, it’s not
possible to have ideal scalability because of communication delays
and idle times. As expected, the explicit formulation of HGKS scales
properly with the increasing number of GPUs. When GPU code using
more than $8$ GPUs, the communication across GPU nodes with RoCE is
required, which accounts for the worse scalability using $16$ GPUs.

Meanwhile, the detailed data for Tesla V100 GPU are given in
Table.\ref{GPU-CPU-TG-D}, in which total computational time and the
percentage of the time for communication are given. Specifically,
the communication time is consist of the time for MPI$\_$ISEND and
MPI$\_$IRECV for the transfer of conservative variables of ghost
cells, MPI$\_$WAITALL for ensuring completion of communications, and
MPI$\_$ALLREDUCE for global time step. The communication across GPU
nodes with RoCE consumes longer time than the communication in
single GPU node with NVLink. The performance of communication across
GPU nodes with InfiniBand will be tested,  which is designed for HPC
centers to improve efficiency among GPU nodes.

\section{Conclusion}
In this paper, to accelerate the computation, HGKS is implemented
with GPU using CUDA on unstructured meshes. For single-GPU
computation, the connectivity of geometric information is generated
for the requirement of data localization and independence. A simple
setting for thread and block is adopted for arbitrary unstructured
meshes, and the executions are implemented automatically by GPU. To
further improve the computational efficiency and enlarge the
computational scale, the HGKS code is also implemented with multiple
GPUs using MPI and CUDA. For multiple-GPU computation, the domain
decomposition and data exchange need to be taken into account.  The
domain is decomposed into subdomains by METIS, and the MPI processes
are created for the control of each process and communication among
GPUs. With reconstruction of connectivity and adding ghost cells,
the main configuration of CUDA for single-GPU can be inherited by
each GPU. For single-GPU implementation using CUDA,  compared with
the CPU code using  2 Intel(R) Xeon(R) Gold 5120 CPUs with OpenMP
directives, 5x speedup is achieved for RTX A5000 and 9x speedup is
achieved for Tesla V100. For multiple-GPU with CUDA and MPI,  the
HGKS is strongly scalable with the increasing number of GPUs. Nearly
linear speedup can be achieved under suitable computational
work-load. Numerical performance shows that the data communication
crossing GPUs through MPI costs the relative little time. The
comparisons between FP32  (single)  precision and FP64  (double)
precision of GPU are also given. The accuracy test is used for
evaluate the effect of precision on computation of HGKS. As
expected, the efficiency can be improved and the memory cost can be
reduced with FP32 precision. Compared with the results of FP64
precision,  the errors of the accuracy increase slightly and the
third-order can be maintained.

\section*{Acknowledgements}
The current research of L. Pan is supported by Beijing Natural
Science Foundation (1232012),  National Natural Science Foundation
of China (11701038)  and the Fundamental Research Funds for the
Central Universities.

\section*{Declaration of competing interest}
The authors declare that they have no known competing financial
interests or personal relationships that could have appeared to
influence the work reported in this paper.

\section*{Data availability}
The data that support the findings of this study are available from
the corresponding author upon reasonable request.

\end{document}